\pdfoutput=1 

\documentclass[applsci,article,accept,pdftex,moreauthors]{Definitions/mdpi} 
\graphicspath{{Figures/}}
\firstpage{1} 
\pubvolume{15}
\issuenum{23}
\articlenumber{12827}
\pubyear{2025}
\copyrightyear{2025}
\externaleditor{Marcin Wojnakowski} 
\datereceived{11 October 2025} 
\daterevised{14 November 2025} 
\dateaccepted{18 November 2025} 
\datepublished{4 December 2025} 
\hreflink{https://{\linebreak} doi.org/10.3390/app152312827} 

\usepackage{comment}
\usepackage{amsmath}
\usepackage{mathtools}
\usepackage{xcolor, soul}
\sethlcolor{yellow} 
\newcommand{\N}{{\mathbb N}}
\let\epsilon\varepsilon
\let\rho\varrho
\let\phi\varphi

\Title{The Basic Reproduction Number for Petri Net Models: A~Next-Generation Matrix Approach}

\TitleCitation{The Basic Reproduction Number for Petri Net Models: A~Next-Generation Matrix Approach}

\Author{Trevor Reckell $^{1,}${*}\orcidA{}, Beckett Sterner $^{2,}${*} and Petar Jevti\'c $^{1}$}

\AuthorNames{Trevor Reckell, Beckett Sterner and Petar Jevti\'c}

\AuthorCitation{Reckell, T.; Sterner, B.; Jevti\'c, P.}

\address{%
$^{1}$ \quad School of Mathematical and Statistical Sciences, Arizona State University, 901 S. Palm Walk, \mbox{Tempe, AZ 85287-1804, USA}; {petar.jevtic@asu.edu}
\\
$^{2}$ \quad School of Life Sciences, Arizona State University, Tempe, AZ 85281-5000, USA}

\corres{Correspondence: trreckel@asu.edu (T.R.); beckett.sterner@asu.edu (B.S.)}

\abstract{The basic reproduction number ($R_0$) is an epidemiological metric that represents the average number of new infections caused by a single infectious individual in a completely susceptible population. The methodology for calculating this metric is well-defined for numerous model types, including, most prominently, Ordinary Differential Equations (ODEs). The basic reproduction number is used in disease modeling to predict the potential of an outbreak and the transmissibility of a disease, as well as by governments to inform public health interventions and resource allocation for controlling the spread of diseases. A Petri Net (PN) is a directed bipartite graph where places, transitions, arcs, and the firing of the arcs determine the dynamic behavior of the system. Petri Net models have been an increasingly used tool within the epidemiology community. However, no generalized method for calculating $R_0$ directly from PN models has been established. Thus, in this paper, we establish a generalized computational framework for calculating $R_0$ directly from Petri Net models. We adapt the next-generation matrix method to be compatible with multiple Petri Net formalisms, including both deterministic Variable Arc Weight Petri Nets (VAPNs) and stochastic continuous-time Petri Nets (SPNs). We demonstrate the method’s versatility on a range of complex epidemiological models, including those with multiple strains, asymptomatic states, and nonlinear dynamics. Crucially, we numerically validate our framework by demonstrating that the analytically derived $R_0$ values are in strong agreement with those estimated from simulation data, thereby confirming the method's accuracy and practical utility.}


\featuredapplication{This work can be utilized by bio-mathematicians, epidemiologists, and other researchers to calculate the basic reproduction number ($R_0$) for the Petri Net system they are modeling, helping them understand the dynamics of an outbreak. }
\begin{document}


\section{Introduction}
\label{intro-sec}
Petri Nets are a promising modeling framework for infectious disease systems that are difficult for standard Ordinary Differential Equation (ODE) models \cite{reckell_numerical_2024,segovia_petri_2025,peleg2005using}. Petri Net models are easily visualized as causal network diagrams, and their modular structure can be scaled to large spatial models \cite{chiaradonna_mpat_2024}. They can also integrate continuous-time models of disease spread with discrete-time events such as policy changes or interventions. Numerous studies have explored the use of Petri Nets to represent and analyze SIR models, highlighting their potential for modular and large-scale modeling of disease dynamics \cite{bahi-jaber_modeling_2003,peng2021modeling,kong_cpn-based_2021, connolly2022epidemic, libkind2022algebraic,segovia_petri_2025}. 

However, Petri Nets have lacked a generalized method for deriving the basic reproduction number, $R_0$. In epidemiology, the basic reproduction number ($R_0$) refers to the expected number of new cases that one infected individual will cause in a population where the entire population is susceptible \cite{dietz1993estimation,diekmann1990definition}. As such, $R_0$ is a fundamental concept in mathematical epidemiology for all communicable diseases. While the Next Generation Matrix (NGM) framework is standard for ODEs, its adaptation to Petri Nets is not yet widely established. Concurrent with our work \cite{reckell_basic_2025}, Segovia has also proposed a formal adaptation of the next-generation matrix method to Petri Nets, focusing on a geometric interpretation of the matrix entries \cite{segovia_petri_2025}. Our paper further advances this emerging topic by establishing a universal, computationally-oriented framework that offers distinct contributions. Compared to the geometric interpretation \cite{segovia_petri_2025}, our framework provides a higher degree of formalization for complex models through a direct method, a wider range of applications, and a more comprehensive verification method. Specifically, our approach is demonstrated on both traditional Stochastic Petri Nets (SPNs) and, crucially, the more flexible Variable Arc Weight Petri Nets (VAPNs). This VAPN compatibility (which can be applied to a wider range of models, such as those in Section \ref{sec:NL}) enables our framework to compactly model complex, nonlinear dynamics that are challenging for path-based methods. Furthermore, we provide, to our knowledge, the first robust numerical verification of an NGM Petri Net framework (Section \ref{sec:Numerical_Verification}), confirming its accuracy and practical utility.  

The inability to directly calculate $R_0$ for Petri Nets limits their analytical utility and poses a significant barrier to the broader adoption of the modular construction of large-scale, mechanistic models of infectious disease spread. An accurate calculation of $R_0$ is therefore crucial for understanding the potential spread of an infectious disease and informing public health interventions \cite{delamater2019complexity}. Analyzing $R_0$ can give insights into the future dynamics of a particular disease based on how its value depends on the specific parameters, subpopulations, and constants that characterize the system at hand. Researchers calculate $R_0$ in models to assess transmission dynamics, compare scenarios, and evaluate the effectiveness of control measures such as vaccination or social behavior \cite{delamater2019complexity}.

Applying the NGM method for calculating $R_0$ within the framework of Petri Nets requires insights  bridging the gap between mechanistic model structure and formal spectral analysis. Mathematical treatments of $R_0$ via the next-generation matrix method have been outlined by Diekmann \cite{diekmann1990definition, diekmann2000mathematical, diekmann2010construction} and Van den Driessche \cite{van2002reproduction}. Although other approaches exist for calculating $R_0$, such as stability analysis of the disease-free equilibrium (DFE)~\cite{van2002reproduction}, the direct definition from survival functions \cite{hethcote2000mathematics}, and estimation from final size relations~\cite{brauer2008mathematical}, the next-generation method is one of the most commonly used for ODE models \cite{delamater2019complexity}.

In this paper, we show how to define and characterize the transition and transmission matrices necessary to calculate $R_0$ in Petri Nets. Among the various formulations of the NGM method, the computational approach introduced by Diekmann et al. \cite{diekmann1990definition, diekmann2000mathematical, diekmann2010construction} is uniquely well-suited for implementation in Petri Nets, as it explicitly constructs the next-generation matrix from compartmental transitions, which naturally aligns with the token flow and transition structure of Petri Nets. This compatibility enables a direct mapping from Petri Net dynamics to the linearized infection and transition matrices necessary for calculating $R_0$. 

\textls[-5]{We demonstrate the applicability and generality of this approach across a range of epidemiological models, including classical structures such as the SIRS, SEEIR models, a compartmental COVID-19 model with various infected compartments, and a model featuring a nonlinear infection term. These examples collectively illustrate how the Diekmann formulation, adapted to Petri Nets, enables the rigorous and automated computation of $R_0$ in Petri Nets, supporting their use across a broader spectrum of infectious disease modeling.}

A premise of this work is that epidemiological models originally formulated as systems of ODEs can be systematically represented as Petri Nets, allowing for a more modular and visual approach to modeling. There are two main methodologies for mapping from PNs to ODEs that facilitate this translation. The first, commonly used for Stochastic Petri Nets (SPNs), maps the system based on mass-action principles. In this approach, outlined by Baez and Biamonte \cite{baez2012quantum}, the state variables (compartments) of the ODE system correspond to the places in the PN, and each term in the differential equations, representing an infection or recovery event, corresponds to a transition governed by a stochastic rate. The second method maps ODEs to a deterministic Variable Arc Weight Petri Net (VAPN). Here, the entire mathematical expression for the rate of change between two compartments can be directly encoded as a functional arc weight on the arc connecting the corresponding places. Both the SPN and VAPN have been numerically verified for the SIRS system \cite{reckell_numerical_2024}. Our framework is intentionally designed to be compatible with both of these ODE to PN mapping formalisms. By demonstrating that our method for determining $R_0$ for Petri Nets yields equivalent $R_0$ values for PN models derived from either method, we ensure its broad applicability and relevance to the wider epidemiological modeling community.

The remainder of the paper is structured as follows. In Section \ref{sec:PN_review}, we introduce the Petri Net formalism and discuss two approaches for using PNs to implement Susceptible--Infectious--Recovered (SIR) models. One approach uses deterministic, discrete-time PNs with variable arc weights, and the other uses stochastic, continuous-time PNs with constant arc weights representing the mass action of interactions among individuals in the compartments. Section \ref{sec:NGM-Pn} presents our NGM framework for Petri Nets. In Section \ref{sec4}, we demonstrate the broad applicability of this framework through a series of case studies on diverse and complex epidemiological models. Finally, a key contribution of this work is presented in Section \ref{sec:Numerical_Verification}, which provides a robust numerical verification of our framework. By simulating the Petri Net models directly and applying established statistical methods to the output data, we demonstrate a strong concordance between our analytically derived $R_0$ and the epidemiologically observed $R_0$ from simulations, thereby validating the accuracy of our approach.

\section{A Brief Overview of Petri Nets}\label{sec:PN_review}

A Petri Net graph, or Petri Net structure, is a weighted bipartite graph \cite{cassandras2008introduction} defined as an $n$-tuple 
    $(P, T,A, w, x)$ where: 
\begin{itemize}
    \item $P$ is the finite set of places (one type of vertex in the graph).
    \item $T$ is the finite set of transitions (the other type of vertex in the graph).
    \item $A$ is the set of arcs (edges) from places to transitions and from transitions to places in the graph $A \subset (P \times T) \cup (T \times P)$. 
    \item $w$ : $A$ → \{1, 2, 3, \ldots\} is the weight function on the arcs.
    \item $x$ is a marking of the set of places $P$; $x_i = [x(p_1), x(p_2), \ldots, x(p_n)$] $\in\N^n$ is the row vector associated with $x$. The term $x_i$ refers to marking of $x$ at state $i$. Note that $x_0$ is the initial state (also known as marking).
\end{itemize}
{Tokens} 
 are assigned to places, with the initial assignment being the initial marking. The number of tokens assigned to a place is an arbitrary non-negative integer, but it does not necessarily have an upper bound. A transition $t_j \in T$ in a Petri Net is said to be enabled if $x(p_i) \geq w(p_i, t_j)$ for all $p_i \in I(t_j)$, where $I(t_j)$ is the set of input arcs from places to $t_j$. This allows us to define the state transition function, $f : \N^n \times T \to \N^n$, of Petri Net $(P, T,A, w, x)$ is defined for transition $t_j \in T$ if and only if $x(p_i) \geq w(p_i, t_j)$ for all $p_i \in I(t_j)$. If $f(x, t_j)$ is defined, then we set $x' = f(x, t_j )$, where $x'(p_i) = x(p_i) - w(p_i, t_j) + w(t_j, p_i), i= 1, \ldots, n$. In simple terms, a transition is enabled if the number of tokens in all places connected to that transition via an incoming arc is greater than or equal to the arc weight for the respective arc connected to the transition.  

An important point to note is that there may be multiple possible PN representations for a given ODE model because there are choices in how to implement the PN model. The most common mapping method involves mapping ODEs to Stochastic Petri Nets (SPNs), as outlined by Baez and Biamonte \cite{baez2012quantum}. An alternative method is to map the ODE to a Variable Arc Weight Petri Net (VAPN) \cite{reckell_numerical_2024}. Variable arc weight Petri Nets can have arc weights that depend on parameters, the number of tokens in places at specific times, and time. VAPNs allow for arc weights to be defined as functions that change over time, which can be linear or nonlinear. This method is deterministic and invertible, allowing for potentially lower computational requirements in parameter fitting and simulations, but it is novel and not yet widely used. 

The inclusion of VAPNs is a key feature of our framework's generality. While SPNs excel at modeling mass-action kinetics, VAPNs can represent complex, nonlinear transmission functions (such as frequency-dependent or behavior-driven rates) as single, functional arc weights. This allows for more elegant and compact models. Specifically, VAPNs can directly model mechanisms that are not based on mass-action principles, a limitation of traditional SPNs. As we will demonstrate, our method of identifying the arcs that contribute to the $\mathcal{F}$ and $\mathcal{V}$ matrices works intrinsically on these functional weights, a task not explicitly addressed by methods tailored to the constant-rate transitions of traditional SPNs.

\section[Next-Generation Matrix for Petri Nets]{{Next-Generation Matrix for Petri Nets}
}\label{sec:NGM-Pn}

We first review the next-generation method for ODEs before demonstrating how it can be applied to PNs. 

\subsection{Next-Generation Matrix for Ordinary Differential Equations} \label{sec:NGM_ODE}
The Next-Generation Matrix (NGM) is a standard and widely-used technique in mathematical epidemiology for calculating $R_0$ from compartmental models. For additional methods for calculating $R_0$, see \cite{martcheva2015introduction}. The method for calculating $R_0$ via the NGM laid out by Diekmann et al.
\cite{diekmann1990definition, diekmann2000mathematical, diekmann2010construction} defines $R_0$ as the spectral radius of the next-generation matrix. Specifically, two matrices (the transmission and transition matrices) are defined to calculate $R_0$, starting with $T$, which is defined as the transmission matrix. The transmission matrix outlines the production of new infections. Where elements within the matrices describe how people move from non-infected compartments to infected compartments. Next, they define $\Sigma$ as the transition matrix, which describes the state changes that occur when individuals move to other infected compartments, die, or become immune. To then find $R_0$ we take the dominant eigenvalue of $(-T\Sigma^{-1})$ yielding $R_0=\varphi(-T\Sigma^{-1})$.

Alternatively, van den Driessche and Watmough \cite{van2002reproduction} outline a similar definition of $R_0$. Let $x=(x_1,…,x_n)^t$, with each $x_i\geq0$, be the number of individuals in each compartment. Let $n$ be the number of compartments, with the first $m$ compartments being the infected compartments. We define $X_s$ as the set of all disease-free states such that $X_s=\{x\geq0|x_i=0,i=1,\ldots,m\}$. Then we define the transmission matrix, matrix $\mathcal{F}_i(x)$, to be the rate of appearance of new infections in compartment $i$. We use two matrices $\mathcal{V}^+_i(x)$ and $\mathcal{V}^-_i(x)$ to describe the rate of transfer of individuals into compartment $i$ by all other means and the rate of transfer of individuals out of compartment $i$, respectively. From $\mathcal{V}^+_i(x)$ and $\mathcal{V}^-_i(x)$, we get the transition matrix, $\mathcal{V}_i=\mathcal{V}^-_i(x)-\mathcal{V}^+_i(x).$ The disease model itself has non-negative initial conditions with the system of {equations following}
\begin{align}
    \dot{x_i} &= f_i(x)= \mathcal{F}_i(x)-\mathcal{V}_i(x), i=1,{\ldots}n.
\end{align}
 While some of these assumptions are trivial when modeled in Petri Nets, others, such as derivations, require closer inspection. For reference, we state the five mathematical assumptions here exactly as they are given in van den Driessche and Watmough \cite{van2002reproduction}, followed by brief comments. We will cover these assumptions for Petri Nets in the following Section \ref{subsec:NGMPN}.

\begin{enumerate}
\item[(1)] If $x\geq0$, then $\mathcal{F}_i,\mathcal{V}^+_i,\mathcal{V}^-_i\geq 0$ for $i=1,…,n$. Each subsequent function $\dot{x_i}$ represents a transfer of individuals, thus they are all non-negative.  
\item[(2)] \textls[-15]{If $x_i=0$ then $\mathcal{V}^-_i=0$. In particular, if $x\in X_s$, then $\mathcal{V}^-_i=0$ for $i=1,…,m$. If a compartment is empty, then there can be no transfer of individuals out of the compartment. }
\item[(3)] $\mathcal{F}_i=0$ if $i>m$. No infectious population enters a non-infectious compartment. 
\item[(4)] If $x\in X_s$ then $\mathcal{F}_i(x)=0$ and $\mathcal{V}^+_i(x)=0$ for $i=1,…,m$. This states that if the population is disease-free, it stays disease-free. Meaning no density-independent infectious population immigrates in. This becomes important when dealing with patch models, as seen in {Section} 
 \ref{sec:Patch}.
\item[(5)] If $\mathcal{F}(x)$ is set to zero, then all eigenvalues of $Df(x_0)$ have negative real parts, where $Df(x_0)$  is the derivative $[\frac{\partial f_i}{\partial x_j}]$ evaluated at the disease free equilibrium (DFE), $x_0$.
\end{enumerate}

These assumptions enable the proof of Lemma 1 from van den Driessche and Watmough \cite{van2002reproduction}. Using Lemma 1, van den Driessche and Watmough subsequently prove that $\varphi(FV^{-1})=R_0$.

\subsection[Constructing the Next-Generation Matrix for Petri Nets]{{Constructing the Next-Generation Matrix for Petri Nets}} \label{subsec:NGMPN}
\begin{figure}[H]
\includegraphics[width=0.9\textwidth]{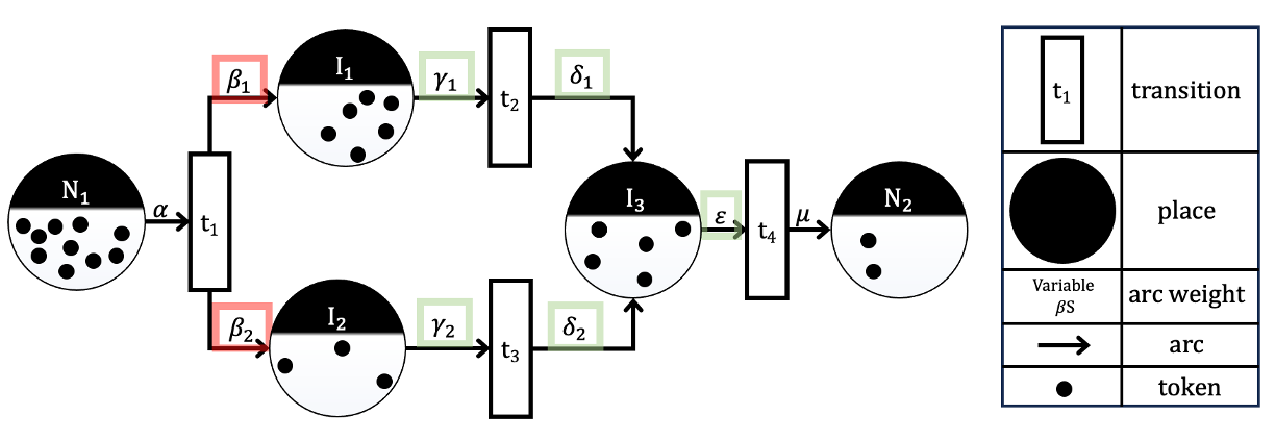}
    \caption{{This} 
 VAPN is used as a guide to help visualize some of the concepts described in the Next-Generation Matrix for Petri Nets formulation. The legend on the right describes the elements within the PN. The places of this VAPN can be labeled $(P_1,P_2,P_3,P_4,P_5)$ or equivalently $(N_1,I_1,I_2,I_3,N_2)$, where the second label gives more description on the infectious state of the place and is used in the~diagram. {Red outline correspond to use in the corresponding $\mathcal{F}$ matrix and green outline correspond to use in the corresponding $\mathcal{V}$ matrix.}}
    \label{fig:R0guide}
\end{figure}

Since many of the papers we draw on for examples below use the van den Driessche and Watmough notation, we will lay out the Next-Generation Matrix for Petri Nets (NGMPN) following similar notation. To help visualize some of the concepts outlined in this process, Figure \ref{fig:R0guide} will be used as a basic VAPN model with non-infected places ($N_1, N_2$) and infected places ($I_1,I_2,I_3$). Examples using SPNs will be provided in the Case Studies section under in the first two model examples. This process of separating disease transmission dynamics from other state transitions within the infected population is conceptually similar to the framework of `Kermack--McKendrick module' proposed by \cite{segovia_petri_2025}. In that work, the Petri Net is partitioned into `susceptible', `infection-process', and `infection' modules, where the $F$ matrix is derived from the infection-process module and the $V$ matrix is derived from the infection module. Our approach implements this separation by identifying specific arcs and transitions corresponding to new infections ($\mathcal{F}_i(x)$) versus those representing recovery, death, or movement between infected states ($\mathcal{V}_i(x)$).
It is important to note that, as with any $R_0$ calculation in any modeling formalism (including ODEs), the modeler must first define the set of infected places or compartments. This is a necessary prerequisite of the epidemiological modeling process itself. Our framework's ``intuitiveness'' lies in how, given this standard partition of places, the $F$ and $V$ matrices can be directly and systematically constructed by identifying the specific arcs that cross this partition, as we will formally define.

In a Petri Net, places are the nodes that hold tokens, representing subpopulations or the state of the system. Then, for the purpose of our framework, it is helpful to think of places in Petri Nets the same way that you might think about compartments in ODEs. Note that in systems biology a ``compartment'' in a Petri Net is a way to group related places and transitions within a larger Petri Net \cite{heiner2009standardised}. Then, we assume we have a Petri Net model for a disease that fits the definition of VAPN or SPN outlined in Section \ref{sec:PN_review}. Assume that each transition can fire once per time step. This assumption does not lead to a loss of generality because the framework's core calculation relies on the rate of change in token flow, which applies to both the functional arc weights in deterministic VAPNs and the continuous-time firing rates in SPNs. The ``time step'' serves as a conceptual unit for a state update rather than a rigid constraint on the model's dynamics. Define the total number of infected places in the Petri Net model as $Q$. Thus, the infected places are defined as $I_1, I_2,\ldots, I_Q$. As a note, we define all places as $P$ and they are labeled $P_{1},\ldots, P_L$, where $L$ is the total number of places and $Q\leq L$. We define $t_{N}$ as a transition with an input arc coming from a non-infected place and an output arc going to an infected place. With $w(t_{N_{j,\ldots,k}}, I_i)$ being the sum of all arcs weights $j,\ldots,k$ from arcs going from non-infected places going to an infected place $I_i$ for VAPN or the net sum of all rate equations $j,\ldots,k$ for arcs from transitions with input arcs from non-infected places and output arcs going to infected places $I_i$ for SPN. Let $x=(x_1,…,x_n)^t$, with each $x_i\geq 0$, be the number of tokens in each place, which is the Petri Net equivalent of population in each compartment in an ODE. Then we define $\mathcal{F}_i(x)$ to be the rate at which tokens move into infected places. For the variable arc weight PN, this rate refers to the arc weights. 

The following provides the clear formal algorithm for distinguishing new infection arcs ($F$) from other state change arcs ($V$). First, let the set of all places be partitioned into infected places, $P_I$, and non-infected places, $P_N$. We then define an infection transition $t \in T$ as any transition that has at least one input arc from a place in $P_N$ and at least one output arc to a place in $P_I$. Let the set of all such infection transitions be $T_{inf}$. From this formal definition, we construct the matrices $F$ and $V$:
\begin{itemize}
\item {For $F$: For any infected place $p_i \in P_I$, we define $\mathcal{F}_i(x)$ as the rate of all new infections entering $p_i$, which is the net inflow from all infection transitions $t_j \in T_{inf}$. The matrix $\mathcal{F}(x)$ is composed of these $\mathcal{F}_i(x)$ terms.} 
\item {For $V$: Conversely, $\mathcal{V}_i(x)$ represents the net rate of all other transfers affecting place $p_i$. This is calculated as the rate of all outflows from $p_i$ (to any place) minus the rate of all inflows to $p_i$ from other infected places $p_k \in P_I$.}
\end{itemize}


The matrices in Equations \eqref{eq:NGMPN_F_Transm} and \eqref{eq:NGMPN_V} provide the operational method for calculating these conceptual rates.
Then $\mathcal{F}_i(x)$ is represented in matrix form by
\begin{align}
\mathcal{F}_{1,\ldots,T}(x) &=\begin{bmatrix}
w(t_{N_{j,\ldots,k}},I_1)& w(t_{N_{j,\ldots,k}},I_1) &  \cdots & w(t_{N_{j,\ldots,k}},I_1)  \\
w(t_{N_{j,\ldots,k}},I_2)& w(t_{N_{j,\ldots,k}},I_2) &  \cdots & w(t_{N_{j,\ldots,k}},I_2) \\
\vdots & \vdots& \ddots& \vdots\\
w(t_{N_{j,\ldots,k}},I_T)& w(t_{N_{j,\ldots,k}},I_T) & \cdots & w(t_{N_{j,\ldots,k}},I_T)) 
\end{bmatrix} .
\label{eq:NGMPN_F_Transm}
\end{align}
Examining Figure \ref{fig:R0guide}, we observe that the parameters of the PN that correspond to the  $\mathcal{F}$ matrix entries are outlined in red. This is continued in all the other Petri Net models in this paper to facilitate easier identification of elements by inspection. For $\mathcal{F}_1$ the element $w(t_{N_{j,\ldots,k}},I_1)=\beta_1$. For $\mathcal{F}_2$ the element $w(t_{N_{j,\ldots,k}},I_2)=\beta_2$. For $\mathcal{F}_3$ no transmission occurs directly into place $I_3$, so the element $w(t_{N_{j,\ldots,k}},I_3)=0$. Thus 
\begin{align*}
\mathcal{F}_{1,2,3}(x) &=\begin{bmatrix}
\beta_1 & \beta_1 & \beta_1\\
\beta_2 & \beta_2 & \beta_2\\
0 & 0 & 0
\end{bmatrix}. 
\end{align*}

We define $w(I_i,t_{j,\ldots,k})$ as the summation of arc weights of output arcs from infected place $I_i$ to transitions $t_{j,\ldots,k}$ for VAPNs and the summation of rate equations of all arcs between infected place $I_i$ to transitions $t_{j,\ldots,k}$. In other words, this represents adding up all arcs leaving the infected place $I_i$. We define $w(t_{I_i}, I_j)$ as the arc weights of arcs from transition to infected places where the arcs leading into the transition come from infected place $I_i$, and this arc weight is equal to the arc weight of the arc from $t_{I_i}$ to $I_j$. Lastly, we define $w(t_{N_m,\ldots,N_l}, I_j)$ as the sum of the arc weights of arcs from transitions $t_{N_m,\ldots,N_l}$ to infected place $ I_j$, where these transitions have no arcs leading into them. We define $V^-_i(x)$ to be the arc weight of tokens leaving infected place $i$ and $V^+_i(x)$ to be the arc weight of tokens going into infected place $i$ outside of coming from non-infected places. In terms of what the terms $\mathcal{V}^-_i(x)$ and $\mathcal{V}^+_i(x)$ represent, it is the rate of people leaving the infected place $i$ and the rate of people going into infected place $i$ from anything other than coming from non-infected places, respectively. These terms combine to yield, $\mathcal{V}_i(x)=\mathcal{V}^-_i(x)-\mathcal{V}^+_i(x)$. For Petri Nets, this can be expressed through the matrix
\begin{equation}
\mathcal{V}_{1,\ldots,T}(x)=\begin{bsmallmatrix}
w(I_1,t_{j,\ldots,k})-w(t_{N_m,\ldots,N_l}, I_1)& -w(t_{I_2}, I_1) &  \cdots & -w(t_{I_T}, I_1)  \\
-w(t_{I_1}, I_2) & w(I_2,t_{j,\ldots,k})-w(t_{N_m,\ldots,N_l}, I_2)  & \cdots & -w(t_{I_T}, I_2) \\
\vdots & \vdots&  \ddots& \vdots\\
-w(t_{I_1}, I_T) & -w(t_{I_2}, I_T) &  \cdots & w(I_T,t_{j,\ldots,k})-w(t_{N_m,\ldots,N_l}, I_T)
\end{bsmallmatrix}. \label{eq:NGMPN_V}
\end{equation}
Examining Figure \ref{fig:R0guide}, we observe that the elements used in $\mathcal{V}$ are outlined in green. This is continued in all the other Petri Net models in this paper to facilitate easier identification of elements that fit in $\mathcal{V}$ by inspection. For $\mathcal{V}_{1}$ the elements $w(I_1,t_{j,\ldots,k})=\gamma_1$,\linebreak  $w(t_{N_m,\ldots,N_l}, I_1)=0$, $-w(t_{I_2}, I_1)=0$, and $-w(t_{I_3}, I_1)=0$. For $\mathcal{V}_{2}$ the elements\linebreak  $-w(t_{I_1}, I_2)=0$, $w(I_2,t_{j,\ldots,k})=\gamma_2$, $w(t_{N_m,\ldots,N_l}, I_2)=0$, and $-w(t_{I_3}, I_2)=0$. For $\mathcal{V}_{3}$ the elements $-w(t_{I_1}, I_T)=-\delta_1$, $-w(t_{I_2}, I_3)=-\delta_2$, $w(I_3,t_{j,\ldots,k})=\epsilon$, and $w(t_{N_m,\ldots,N_l}, I_3)=0$.  Thus, 
\begin{align*}
\mathcal{V}_{1,2,3}(x) &= \begin{bmatrix}
\gamma_1 & 0 & 0\\
0 & \gamma_2 & 0\\
-\delta_1 & -\delta_2 & \epsilon
\end{bmatrix}.
\end{align*}

The disease free equilibrium (DFE) for a system with $L$ places is defined as ($P_1^*,\ldots, P_L^*$), the state at which the system remains free of the disease. For the Petri Nets, the DFE can be found by first setting the token level of all infected places to zero. Then, the arc weights of arcs into and out of a given non-infected place $P_i$ are summed, with the arc weights of arcs going out of $P_i$ being multiplied by $-$1. The resulting expression from these steps is then set to zero. From here, a system of equations is formed that can be solved for $P_1^*,\ldots, P_L^*$. Additionally, the total population, usually denoted as $N$, is defined as $\mathcal{N}+\epsilon$, where $\mathcal{N}$ is the summation of tokens in all places, and $\epsilon$ is the error from rounding in a discrete token level Petri Net. In a continuous token level Petri Net, $\epsilon=0$. The difference in total population between $N$ and $\mathcal{N}+\epsilon$ becomes relevant when performing actual simulations of models, rather than when finding $R_0$ using the NGMPN method. From here, a system of equations is formed where the DFE token levels labeled for each non-infected place can be solved for. 

Once these $\mathcal{F}_i(x)$, $\mathcal{V}_i(x)$, and DFE values are obtained, they are then converted into the respective $F$ (the transmission matrix) and $V$ (the transition matrix) Jacobian matrix by finding the rate of change with respect to $I_1$ for column one, $I_2$ for column two, $\cdots$, $I_T$ for column  $T$, at the DFE, for the respective arc weights defined as functions. 
Note that Petri Nets can be in discrete or continuous time. Thus, the rate of change would be the difference quotient or derivative, respectively. 

Then, using the Petri Net-based $F$ and $V$ values, the same Next-Generation Matrix steps of finding the inverse of $V$ and finding the dominant eigenvalues of $FV^{-1}$ can be done to find the Petri Net $R_0$. If the values of $F, V$, and the DFE of two systems are equal, then the resulting $R_0$ values are equal. However, the way parameters interact within an ODE model and a PN model differ. We will demonstrate the process in detail, providing various examples with different functionalities.

Addressing the assumptions (A1)--(A5) made in van den Driessche and Watmough~\cite{van2002reproduction}, we interpret them for Petri Nets. 
\begin{enumerate}[leftmargin=3em,labelsep=3.5mm]
\item[(A1)]  If a place $P_i$ contains a non-negative number of tokens $t_i\geq 0$, then the number of tokens that can be added to or removed from any place by any transition is non-negative.
Specifically, if a transition $t_j$ is enabled (based on its input places and corresponding markings), the number of tokens removed from the input places and added to the output places must be non-negative. By definition of Petri Nets, the assumption is inherently preserved across all Petri Net structures. This is due to the stipulation that the number of tokens in each place must remain non-negative at all times. Additionally, transitions in a Petri Net operate in a manner that ensures the non-negativity of tokens, effectively maintaining this essential characteristic.   \label{A1}

\item[(A2)] If the marking of a place $P_i$ is zero (i.e., $t_i=0$), then no transition can consume tokens from place $P_i$. Specifically, if a place is empty (i.e., no tokens are present in the corresponding place), no transition can have it as an input. This ensures that no ``removal'' of individuals (tokens) from an empty place occurs. Again, the definition of Petri Net preserves this assumption for all Petri Nets, as a transition can not fire if an input place has less tokens than the corresponding arc weight. If a place has no tokens (individuals), no transition can remove tokens from that place.

\item[(A3)]If place $P_i$ corresponds to a non-infectious place, then no transition can place tokens into this place from a transition that represents infection occurring or remaining. No transition should be able to fire that would place tokens into non-infectious places from a transition that represents infection occurring or remaining. 

\item[(A4)] If the marking of a place $P_i$ indicates a disease-free population (i.e., there are no tokens representing infected individuals in this marking), then no transition can fire that would result in tokens (representing infected individuals) being moved to an infected place. If a place has a disease-free marking, no transitions can place tokens in it that correspond to infected individuals (i.e., no ``infection'' transition should have this place as an output if it is currently disease-free).


\item[(A5)]   If $F(x)$ is set to $0$, then all eigenvalues of the Jacobian of the infected subsystem at the DFE have negative real parts. In terms of the Petri Nets, this means that when no new infections occur, that means all transitions transferring tokens from non-infected to infected places are disabled, the marking of every infected place will asymptotically decay to zero. Formally, let the state vector be partitioned as $x = (x_I, x_S)$, where $x_I$ are the infected places and $x_S$ are the non-infected places. Define the infected subsystem as $\dot{x}_I = F_I(x_I, x_S) - V_I(x_I, x_S)$, where $F_I$ collects all token inflows into infected places from non-infected places and $V_I$ collects all other token flows (recoveries, deaths, or transfers among infected places). Then, at the DFE $x_0$, with $F_I(x_0) = 0$, the Jacobian matrix 
\[
A = D_{x_I}(F_I - V_I)\big|_{x = x_0} = -D_{x_I}V_I(x_0)
\]
\textls[-15]{must be Hurwitz, which means all eigenvalues have strictly negative real parts. Note that $D$ is the derivative operator. Equivalently, $D_{x_I}V_I(x_0)$ should be a nonsingular $M$-matrix, ensuring that small perturbations in the infected places decay over time. In terms of a biological application, this assumption guarantees that in the absence of new infections, any initial infection will die out, preserving the disease-free equilibrium. \label{A5}}
\end{enumerate} 

\section{Case Studies}\label{sec4}
We have selected a range of ODE and corresponding PN example models to demonstrate how the Next-Generation Matrix for Petri Nets approach to determining $R_0$ works. Each model either adds additional compartments/places or introduces new dynamics, providing a range of examples for people to draw on if they wish to implement the NGMPN in a new model.

\subsection{{SIRS}
} \label{sec:SIRS}

We start with a simple example of an $SIRS$ system, a slight variation in the Kermack--McKendrick model \cite{kermack1927contribution}. This shows a straightforward application of the NGMPN process for a foundational epidemiological model. The compartments $S$, $I$, and $R$ represent susceptible, infected, and recovered, respectively. The parameters of this system $\beta, \delta,$ and $\gamma$ represent the rate of re-susceptibility, rate of infection, and rate of recovery, respectively. The SIRS system is described by the following system of equations:
\begin{align}
\frac{dS}{dt} &= -\beta SI+\delta R,\label{eq:S_SIRS} \\
\frac{dI}{dt} &= \beta S I- \gamma I, \label{eq:I_SIRS} \\
\frac{dR}{dt} &= \gamma I-\delta R. \label{eq:R_SIRS} 
\end{align}

\begin{figure}[H]    
    \subfloat[\centering]{\includegraphics[width=0.49\textwidth]{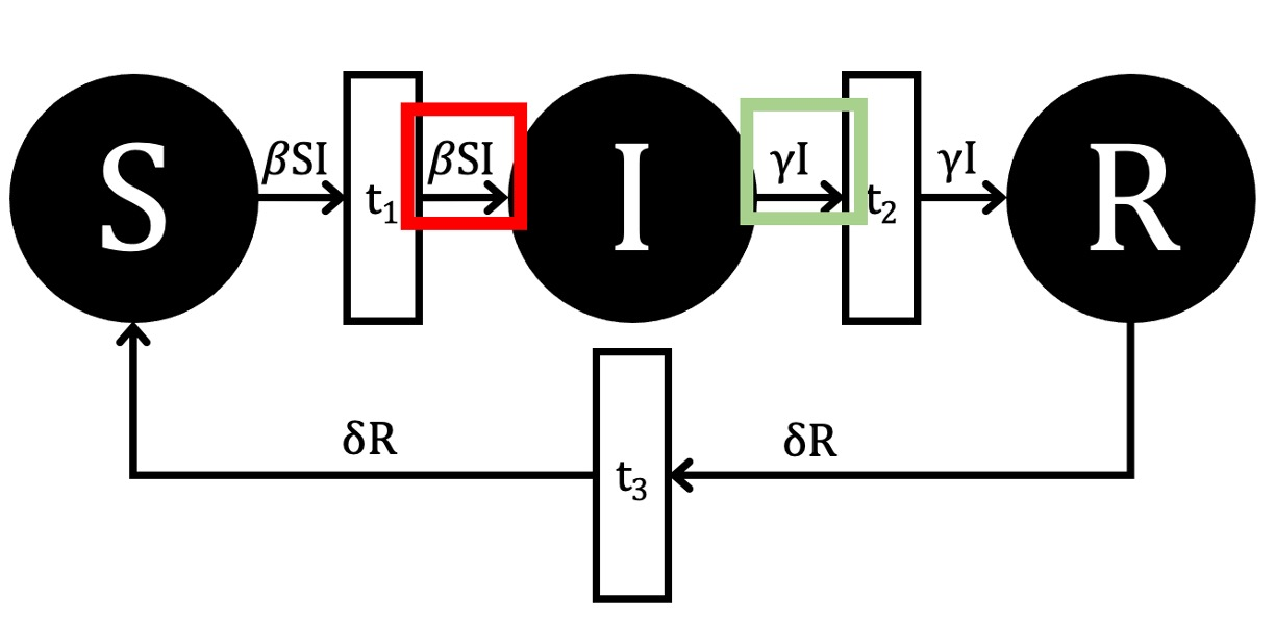}
    }
    \hfill
    \subfloat[\centering]{\includegraphics[width=0.49\textwidth]{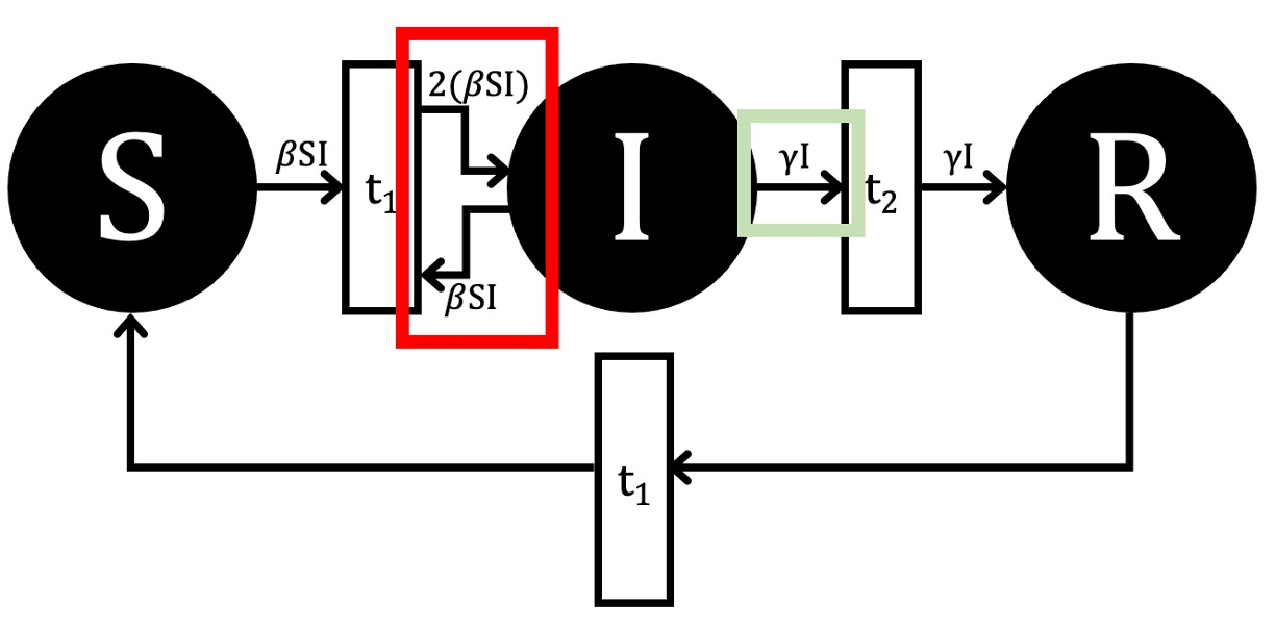}
    }
    \caption{{Two Petri} 
 Net implementations of the SIRS ODE model given in Equations \eqref{eq:S_SIRS}--\eqref{eq:R_SIRS}. Option (\textbf{a}) uses a discrete-time, variable arc weight PN. Option (\textbf{b}) uses a continuous-time, stochastic Markov Chain PN. {The arcs and arc weights used in forming matrix $F_i(x)$ are outlined and highlighted in red, and the arcs and arc weights used in forming matrix $V_i(x)$ are outlined and highlighted in green.}}
    \label{fig:SIRS}
\end{figure}

The ODE system can be mapped into two Petri Net models: one with variable arc weights and minimal arcs, yielding Figure \ref{fig:SIRS}a, and another using the methods outlined by Baez and Biamonte \cite{baez2012quantum}, yielding Figure \ref{fig:SIRS}b.

\vspace{-17pt}

The first step in finding the $R_0$ of the PN systems in Figure \ref{fig:SIRS} is to obtain the Disease Free Equilibrium (DFE) by setting the token levels to zero in the infected places, setting the net arc weight of non-infected places to 0, and solving for the token values in non-infected places $(S^*, I^*, R^*)=(S^*, 0, 0)$, where $S^*=N$. Where $N=S+I+R$ or the sum of the entire population. Then the transmission matrix starts with finding $\mathcal{F}_i(x)$, which can be found by looking at the arc weights of arcs going into the infected compartments as shown in Equation \eqref{eq:NGMPN_F_Transm}. These arcs are outlined in red in Figure \ref{fig:SIRS}. Note that $i$ refers to the corresponding infected compartment, which we have numbered. For this example, the $I$ compartment is 1. From here we find $w(t_{N_{j,\ldots,k}},I_1)=\beta SI$ for the VAPN in Figure \ref{fig:SIRS}a and $w(t_{N_{j,\ldots,k}},I_1)=2\beta SI-\beta SI$ for the SPN in Figure \ref{fig:SIRS}b. Thus,
 \begin{align*}
\mathcal{F}_1(x) &=\begin{bmatrix}
\beta SI
\end{bmatrix},
\end{align*}
for the Petri Net layout as depicted in Figure \ref{fig:SIRS}a, and
 \begin{align*}
\mathcal{F}_1(x) &=\begin{bmatrix}
2\beta SI-\beta SI
\end{bmatrix}. 
\end{align*}
for the Petri Net layout as depicted in Figure \ref{fig:SIRS}b.
From this, we find the rate of change with respect to the infected place $I$, and we evaluate it at the DFE to have
\begin{align*}
F &=\begin{bmatrix}
\beta N
\end{bmatrix}.
\end{align*}
for both layouts.
For the transition matrix, we start with $\mathcal{V}$, and we look at arc weights of arcs leaving the infected places as shown in Equation \eqref{eq:NGMPN_V}. These arcs are outlined in green in Figure \ref{fig:SIRS}. The terms $w(I_1,t_{j,\ldots,k})=\gamma I$ and $w(t_{N_m,\ldots,N_l}=0$ for both the VAPN in Figure \ref{fig:SIRS}a and the SPN in Figure \ref{fig:SIRS}b. Thus yielding for both layouts
\begin{align*}
\mathcal{V}_1(x) &= \begin{bmatrix}
\gamma I
\end{bmatrix}.
\end{align*}
We then find the rate of change with respect to the infected place $I$ and evaluate at the DFE. From this, we find the inverse:
\begin{align*}
V &= \begin{bmatrix}
\gamma 
\end{bmatrix}, \\
V^{-1} &= \begin{bmatrix}
\frac{1}{\gamma} 
\end{bmatrix}.
\end{align*}
Now that we have the $F$ and $V$ matrices, we can find $R_0$ by multiplying them and then determining the dominant eigenvalue we get,
\begin{align*}
R_0 &= \varphi\begin{bmatrix}
F V^{-1} 
\end{bmatrix} = \frac{\beta N}{\gamma}.
\end{align*}

This $R_0$ value of both layouts of  $\frac{\beta N}{\gamma}$ is analytically equivalent to the ODE $R_0$ value.

\subsection{SEIR} \label{sec:SEIR}
This second example introduces an exposed compartment for diseases with delayed onset. Cooke laid out a version of the  model Cooke in 1966 \cite{cooke1966functional}. The additional compartment $E$ represents exposed individuals, usually meaning a person is infected with a disease but cannot infect others yet. The parameters $\Pi, \beta, \mu, \eta, $ and $\alpha$ represent the birth/immigration rate, the infection rate, the death/emigration rate, the onset of infectiousness/symptoms rate, and the recovery rate, respectively. The SEIR system is described by the following system of equations:
\begin{align}
\frac{dS}{dt} &= \Pi-\beta SI -\mu S,\label{eq:S_SEIR} \\
\frac{dE}{dt} &= \beta SI -\eta E-\mu E,\label{eq:E_SEIR} \\
\frac{dI}{dt} &= \eta E-\alpha I -\mu I,\label{eq:I_SEIR} \\
\frac{dR}{dt} &= \alpha I-\mu R.\label{eq:R_SEIR} 
\end{align}

\begin{figure}[H]
    
    \subfloat[\centering]{\includegraphics[width=0.49\textwidth]{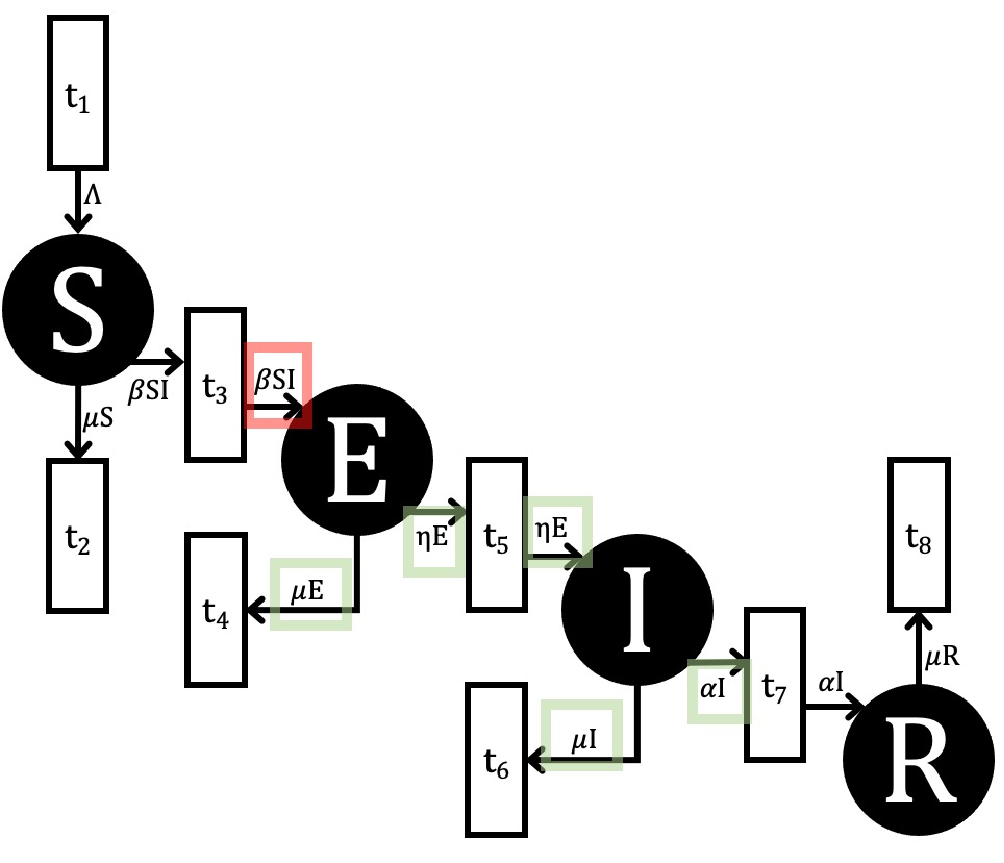}
    }
    \hfill
    \subfloat[\centering]{\includegraphics[width=0.49\textwidth]{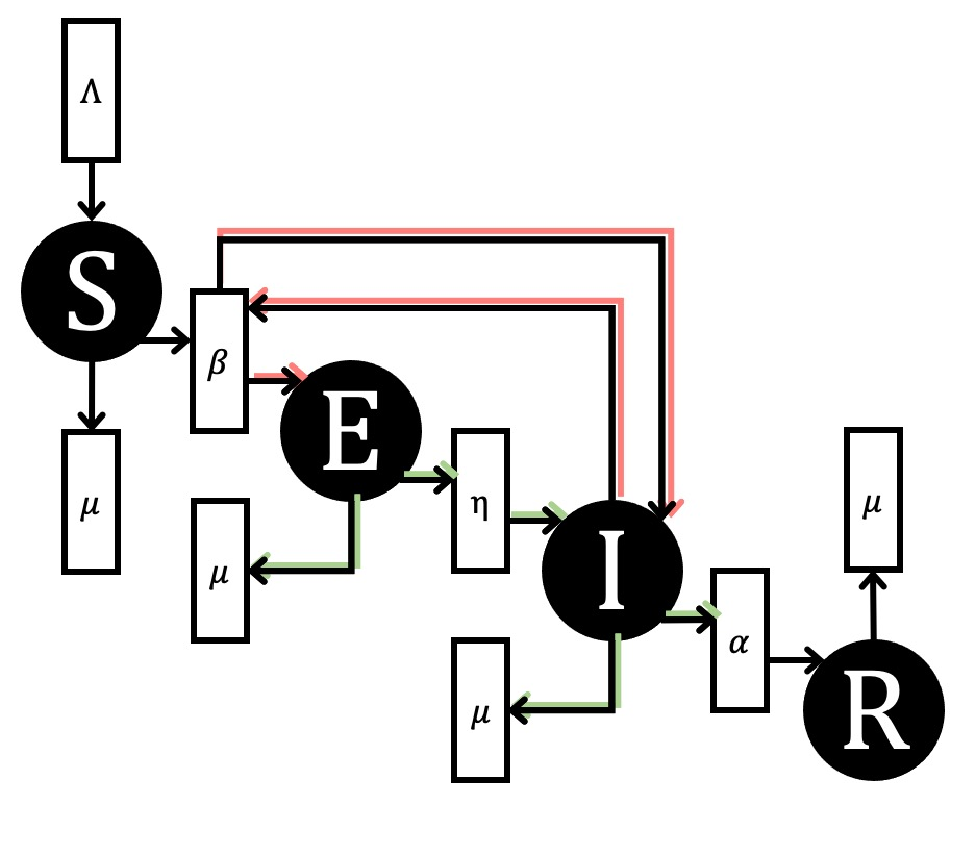}
    }
    \caption{Two Petri Net implementations of the SEIR ODE model given in Equations \eqref{eq:S_SEIR}--\eqref{eq:R_SEIR}. Option (\textbf{a}) uses a discrete-time, variable arc weight PN. Option (\textbf{b}) uses a continuous-time, stochastic Markov Chain PN. {The arcs and arc weights used in forming matrix $F_i(x)$ are outlined and highlighted in red, and the arcs and arc weights used in forming matrix $V_i(x)$ are outlined and highlighted in green.}}
    \label{fig:SEIR}
\end{figure}

The Equations \eqref{eq:S_SEIR}--\eqref{eq:R_SEIR} are mapped to a VAPN (Figure \ref{fig:SEIR}a) and stochastic Petri Net (Figure \ref{fig:SEIR}b). 

\vspace{-12pt}


To find the $R_0$ of the PN systems in Figure \ref{fig:SEIR} we follow the same procedure outlined in Section \ref{sec:SIRS}, we first obtain the DFE. By setting the infected tokens to zero and solving for the non-infected steady state, we obtain $(S^*, E^*, I^*, R^*)=(\frac{\Pi}{\mu}, 0, 0, 0, 0)$. Then, the matrix $F_i(x)$ can be found by examining the arc weights of arcs entering the infected compartments, outlined in red in Figure \ref{fig:SEIR}, as shown in Equation \eqref{eq:NGMPN_F_Transm}. Note that $i$ refers to the corresponding infected compartment and the corresponding row in the matrix, which we have numbered. For this example, the $E$ compartment is 1 and the $I$ compartment is 2. The terms for the VAPN in Figure \ref{fig:SEIR}a within $\mathcal{F}_{1}$ are found by inspection to be $w(t_{N_{j,\ldots,k}},I_1)=\beta SI$ and within $\mathcal{F}_{2}$ to be $w(t_{N_{j,\ldots,k}},I_1)=0$. Similarly the terms for the SPN in Figure \ref{fig:SEIR}b within $\mathcal{F}_{1}$ are found by inspection to be $w(t_{N_{j,\ldots,k}},I_1)=\beta SI$ and within $\mathcal{F}_{2}$ to be $w(t_{N_{j,\ldots,k}},I_1)=\beta SI-\beta SI=0$. Thus,
 \begin{align*}
\mathcal{F}_{1,2}(x) &=\begin{bmatrix}
\beta SI & \beta SI\\
0& 0
\end{bmatrix}, 
\end{align*}
for VAPN as shown in Figure \ref{fig:SEIR}a, and 
 \begin{align*}
\mathcal{F}_{1,2}(x) &=\begin{bmatrix}
\beta SI & \beta SI\\
\beta SI -\beta SI & \beta SI -\beta SI 
\end{bmatrix}, 
\end{align*}
for the stochastic PN as shown in Figure \ref{fig:SEIR}b. 

As $E$ is our first infected compartment and $I$ is the second, we find the rate of change with respect to $E$ in the first column and $I$ in the second column, and evaluate at the DFE. This gives the same result for both Petri Net~options:
\begin{align*}
F &=\begin{bmatrix}
0 & \beta \frac{\Pi}{\mu} \\
0 & 0 \\
\end{bmatrix}.
\end{align*}

For $\mathcal{V}$, we look at the weights of arcs leaving the infected places and transitioning from one infected place to another, as outlined in green in Figure \ref{fig:SEIR}. Then for $\mathcal{V}_{1}$ the terms $w(I_1,t_{j,\ldots,k})=\eta E+\mu E$, $w(t_{N_m,\ldots,N_l}=0, I_1)$, and $-w(t_{I_2}, I_1)=0$. For $\mathcal{V}_{2}$ the terms $-w(t_{I_1}, I_2)=-\eta E$,  $w(I_2,t_{j,\ldots,k})=\alpha I+\mu I$, and $w(t_{N_m,\ldots,N_l}, I_2)=0 $. Both $\mathcal{V}_{1}$ and $\mathcal{V}_{2}$ are the same for the VAPN and SPN in Figures \ref{fig:SEIR}a and \ref{fig:SEIR}b respectively. Thus,
\begin{align*}
\mathcal{V}_{1,2}(x) &= \begin{bmatrix}
\eta E + \mu E & 0 \\
-\eta E & \alpha I+\mu I 
\end{bmatrix}.
\end{align*}
Next, we find the rate of change with respect to the infected place $E$ in the first column and $I$ in the second column, and we evaluate it at the DFE. From this, we find the inverse\linebreak  to give
\begin{align*}
V &= \begin{bmatrix}
\eta  + \mu  & 0 \\
-\eta & \alpha +\mu 
\end{bmatrix} ,\\
V^{-1} &= \begin{bmatrix}
\frac{1}{\eta +\mu} & 0 \\
\frac{\alpha}{(\eta+\mu)(\alpha+\mu)} & \frac{1}{\alpha+\mu} 
\end{bmatrix}.
\end{align*}
Now that we have the $F$ and $V$ matrices, we can find $R_0$ by multiplying them and then determining the dominant eigenvalue. Thus,
\begin{align*}
R_0 &= \varphi\begin{bmatrix}
F V^{-1} 
\end{bmatrix}
=\frac{\beta \Pi \eta}{\mu(\alpha+\mu)(\eta+\mu)}.
\end{align*}

This is analytically the same $R_0$ expression that can be derived for the ODE using the next-generation matrix method. This demonstrates a key strength of our framework, which is its ability to extract the fundamental epidemiological dynamics relevant to $R_0$ irrespective of the specific Petri Net formalism used for implementation, whether deterministic VAPNs or stochastic SPNs.

\subsection{SEEIR} \label{sec:SEEIR}
This third example introduces parallel paths of infected compartments. The model and subsequent $R_0$ value were outlined by Diekmann \cite{diekmann2010construction}. The parameters $\mu, \beta, p, v_1, v_2, $ and $ \gamma$ represent birth/death rate, transmission rate, rate of exposure to infection, rate of becoming infectious for $E_1$, rate of becoming infectious for $E_2$, and rate of recovery, respectively. The SEEIR system is described by the following system of equations:
\begin{align}
\frac{dS}{dt} &= \mu N -\beta \frac{SI}{N}-\mu S,\label{eq:S_EEIR} \\
\frac{dE_1}{dt} &= p\beta \frac{SI}{N}-v_1 E_1 -\mu E_1,\label{eq:E1_EEIR} \\
\frac{dE_2}{dt} &= (1-p)\beta \frac{SI}{N}-v_2 E_2 -\mu E_2,\label{eq:E2_EEIR} \\
\frac{dI}{dt} &= v_1 E_1 + v_2 E_2 -\gamma I -\mu I,\label{eq:I_EEIR} \\
\frac{dR}{dt} &= \gamma I-\mu R. \label{eq:R_EEIR} 
\end{align}

\begin{figure}[H]
    \includegraphics[width=0.7\textwidth]{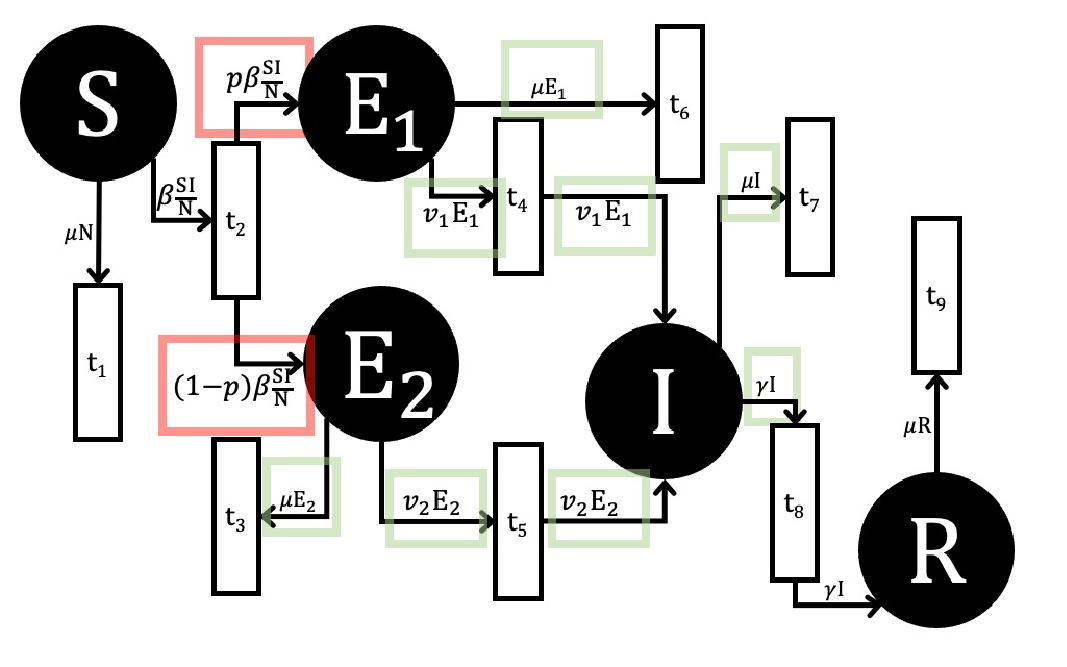}
    \caption{This SEEIR Petri Net model is mapped directly from Equations \eqref{eq:S_EEIR}--\eqref{eq:R_EEIR}. {The arc weights used in forming matrix $F_i(x)$ are outlined in red, and the arc weights used in forming matrix $V_i(x)$ are outlined in green.}}
    \label{fig:SEEIR_PN}
\end{figure}

The Equations \eqref{eq:S_EEIR}--\eqref{eq:R_EEIR} are mapped to a variable arc weight Petri Net to yield \mbox{Figure \ref{fig:SEEIR_PN}}. We will analyze only the Variable Arc Weight Petri Net (VAPN) implementation going forward, for brevity, as it is the more novel approach. The corresponding SPN models would be analyzed using the same principles, relying on the net sum of continuous-time rate equations for each transition, as demonstrated for the SEIR model, yielding identical final $F$ and $V$ matrices at the DFE.

\vspace{-3pt}

The first step in finding the $R_0$ of the PN system in Figure \ref{fig:SEEIR_PN} is to obtain the DFE by setting the token levels of the infected places to zero, setting the net arc weight of non-infected places to 0, and solving for the token values at the non-infected places $(S^*, E_1^*, E_2^*, I^*, R^*)=(S^*, 0, 0, 0, 0)$. Note that since $N=S+E_1+E_2+I+R$, $S^*=N$, where $N$ is the total population. Then, the matrix $F_i(x)$ can be found by examining the arc weights of arcs entering the infected compartments of $E_1, E_2,$ and $I$, as shown in Equation~\eqref{eq:NGMPN_F_Transm}. Note that $i$ refers to the corresponding infected compartment and the corresponding row in the matrix, which we have numbered. For this example, the $E_1$ compartment is 1, the $E_2$ compartment is 2, and the $I$ compartment is 3.

 \begin{align*}
\mathcal{F}_{1,2,3}(x) &=\begin{bmatrix}
p\beta \frac{SI}{N} & p\beta \frac{SI}{N} & p\beta \frac{SI}{N}\\
(1-p)\beta \frac{SI}{N} & (1-p)\beta \frac{SI}{N} & (1-p)\beta \frac{SI}{N}\\
0 & 0 & 0
\end{bmatrix}. 
\end{align*}
From this, we find the rates of change with respect to $E_1$ in the first column, $E_2$ in the second column, $I$ in the third column, and we evaluate them at the DFE to give
\begin{align*}
F &=\begin{bmatrix}
0 & 0 & p\beta\\
0 & 0 & (1-p)\beta\\
0 & 0 & 0
\end{bmatrix}.
\end{align*}
{For $\mathcal{V}$, we analyze the Petri net to determine arcs leaving infected places as shown in Equation} 
 \eqref{eq:NGMPN_V}. This yields
\begin{align*}
\mathcal{V}_{1,2,3}(x) &= \begin{bmatrix}
v_1 E_1 + \mu E_1 & 0 & 0\\
0 & v_2 E_2 + \mu E_2 & 0\\
-v_1E_1 & -v_2E_2 & \gamma I + \mu I
\end{bmatrix}.
\end{align*}
{We then find the rate of change with respect to the infected place $I$, and evaluate at the DFE. From this we find the inverse to give}
\begin{align*}
V &= \begin{bmatrix}
v_1+\mu & 0 & 0\\
0 & v_2+\mu & 0\\
-v_1 & -v_2 & \gamma+\mu 
\end{bmatrix}, \\
V^{-1} &= \begin{bmatrix}
\frac{1}{v_1+\mu} & 0 & 0\\
0 & \frac{1}{v_2+\mu} & 0\\
\frac{v_1}{(v_1+\mu)(\gamma+\mu)} & \frac{v_2}{(v_2+\mu)(\gamma+\mu)} & \frac{1}{\gamma+\mu}
\end{bmatrix}.
\end{align*}
{Now that we have $F$ and $V$ matrices, we can find $R_0$ by multiplying and finding the dominant eigenvalue. Thus}
\begin{align*}
R_0 &= \varphi\begin{bmatrix}
F V^{-1} 
\end{bmatrix} = \left(\frac{pv_1}{v_1+\mu} + \frac{(1-p)v_2}{v_2+\mu}\right) \frac{\beta}{\gamma+\mu}.
\end{align*}

This is analytically the same $R_0$ expression found in Equation 2.10 of \cite{diekmann2010construction}.

\subsection[SVEIR]{SVEIR}  \label{sec:SVEIR}
\textls[-5]{While $R_0$ is the average number of secondary infections produced by a single infectious individual introduced into a completely susceptible population \cite{diekmann2010construction}, there is also the notion of effective reproduction number of $R_e$ which is the average number of secondary infections produced by a single infectious individual at a specific point in time ($t$) \cite{cori2013new}. With this in mind, we move to an extension of the SEIR model by Bugalia et al. with a compartment for Vaccinated individuals $V$ and $N$ is the sum of all compartments \cite{bugalia2023estimating}. With the compartment $V$ we can look at the effective reproduction number when vaccination and non-pharmaceutical interventions are occurring. The parameters $\beta_t, \xi_t, \gamma, $ and $\varphi$ represent the contact rate of individuals (time dependent), vaccination rate (time dependent), the disease re-susceptibility rate, and the vaccine re-susceptibility rate, respectively. 
The parameters $\alpha,$ and $\delta$ represent the recovery rate of infected individuals and disease death rate, respectively. The vaccine efficacy is described as $1-\sigma$ where vaccinated people are infected at the rate $\sigma \beta_t$.
The SVEIR system is described by the following system of equations:}

\begin{align}
\frac{dS}{dt} &= -\frac{S}{N}\beta_t I+\gamma R -\xi_t S+ \varphi V \\ \label{eq:S_SVEIR} 
\frac{dV}{dt} &=  -\frac{V}{N}\sigma \beta_t I + \xi_t S- \varphi V \\
\frac{dE}{dt} &= \beta_t I (\frac{S}{N}+\sigma \frac{V}{N})-\alpha E\\
\frac{dI}{dt} &= \alpha E- (\alpha+\delta) I\\
\frac{dR}{dt} &= \alpha E -\gamma R \label{eq:R_SVEIR} 
\end{align}

\begin{figure}[H]
    \includegraphics[width=0.7\textwidth]{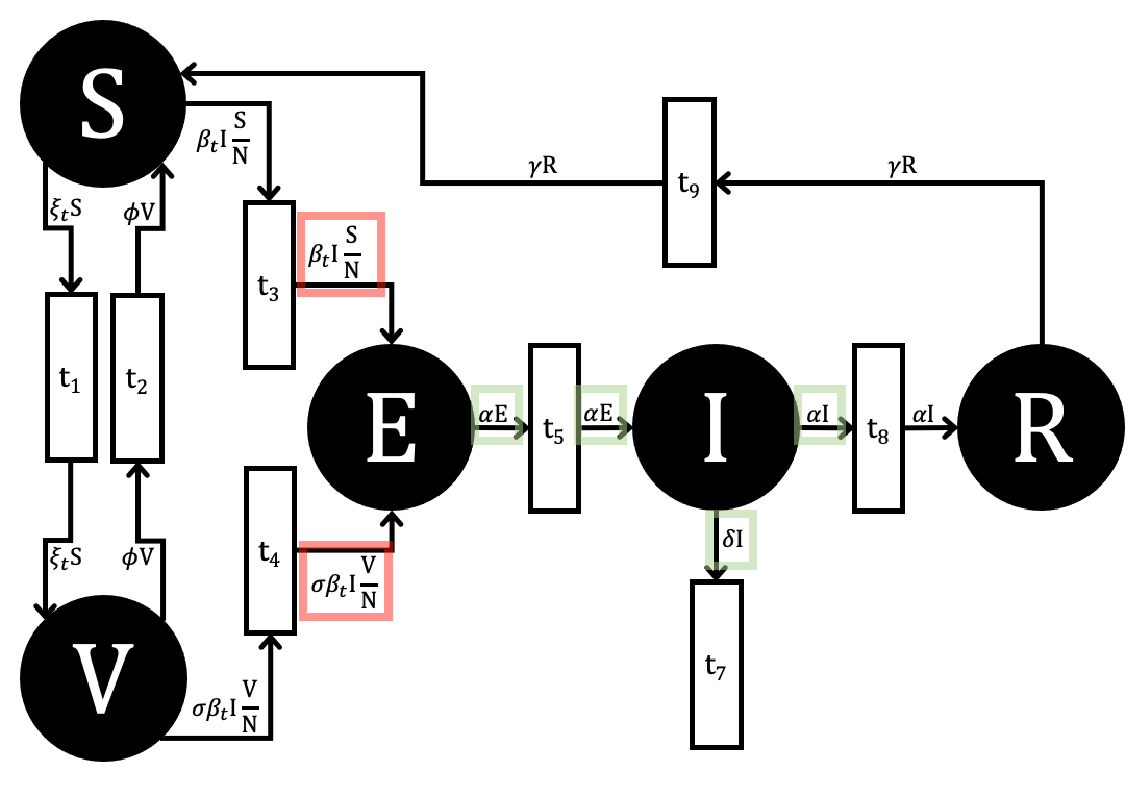}
    \caption{A variable arc weight SVEIR Petri Net mapped from Equations \eqref{eq:S_SVEIR}--\eqref{eq:R_SVEIR}. {The arc weights used in forming matrix $F_i(x)$ are outlined in red, and the arc weights used in forming matrix $V_i(x)$ are outlined in green.}}
    \label{fig:SVEIR_VAW}
\end{figure}

The Equations \eqref{eq:S_SVEIR}--\eqref{eq:R_SVEIR} are mapped to a VAPN in Figure \ref{fig:SVEIR_VAW}. 

\vspace{-7pt}

The steps for finding the effective reproduction number $R_e$ of the PN systems remains the same as the reproduction number for ODEs if Lemma 1 from van den Driessche and Watmough \cite{van2002reproduction} holds true. If the assumptions we laid out for Petri Nets {(A1--A5)} 
 also hold true for our Petri Net model, we then find the effective reproduction number for Petri Nets in the same method as we would for the reproduction number. Thus, we first obtain the DFE by setting the infected places token levels to zero, setting the net arc weight of non-infected places to 0, and solving for the token values at the non-infected places $(S^*, V^*, E^*, I^*, R^*)=(N\frac{\varphi}{\varphi+\xi_t}, N\frac{\xi_t}{\varphi+\xi_t}, 0, 0, 0)$.  Note that since $N=S+V+E+I+R$, $S^*+V^*=N$, where $N$ is the total population. Then, the matrix $F_i(x)$ can be found by examining the arc weights of arcs entering the infected compartments of $E$ and $I$, as shown in Equation \eqref{eq:NGMPN_F_Transm}.
Note that $i$ refers to the corresponding infected compartment and the corresponding row in the matrix, which we have numbered. For this example, the $E$ compartment is 1 and the $I$ compartment is 2. Thus,
 \begin{align*}
\mathcal{F}_{1,2}(x) &=\begin{bmatrix}
\beta_t I(\frac{S}{N}+\sigma\frac{V}{N}) & \beta_t I(\frac{S}{N}+\sigma\frac{V}{N})\\
0& 0
\end{bmatrix}, 
\end{align*}
for VAPN as shown in Figure \ref{fig:SVEIR_VAW}.

We then find the rate of change with respect to $E$ in the first column and $I$ in the second column, and evaluate at the DFE, yielding
\begin{align*}
F &=\begin{bmatrix}
0 & \beta_t (\frac{\varphi+\sigma \xi_t}{\varphi+\xi_t}) \\
0 & 0 \\
\end{bmatrix}.
\end{align*}

For $\mathcal{V}$, we find arcs coming out of infected places, outlined in green to help. This yields
\begin{align*}
\mathcal{V}_{1,2}(x) &= \begin{bmatrix}
\alpha E & 0 \\
-\alpha E & (\alpha+\delta) I
\end{bmatrix}.
\end{align*}

We then find the rate of change with respect to the infected place $I$, and evaluate at the DFE. From this we find the inverse to give
\begin{align*}
V &= \begin{bmatrix}
\alpha & 0 \\
-\alpha & \alpha+\delta
\end{bmatrix}, \\
V^{-1} &= \begin{bmatrix}
\frac{1}{\alpha} & 0 \\
\frac{1}{\alpha+\delta} & \frac{1}{\alpha+\delta}
\end{bmatrix}.
\end{align*}
With the $F$ and $V$ matrices in hand, we can find $R_0$. Thus,
\begin{align*}
R_0 &= \varphi\begin{bmatrix}
F V^{-1} 
\end{bmatrix} = \varphi\begin{bmatrix}
\frac{\beta_t}{\alpha+\delta}(\frac{\varphi+\sigma \xi_t}{\varphi+ \xi_t}) & \frac{\beta_t}{\alpha+\delta}(\frac{\varphi+\sigma \xi_t}{\varphi+ \xi_t}) \\
-\alpha & \alpha+\delta
\end{bmatrix} = \frac{\beta_t}{\alpha+\delta}(\frac{\varphi+\sigma \xi_t}{\varphi+ \xi_t})
\end{align*}

This is analytically the same $R_e$ expression found in Bugalia et al. \cite{bugalia2023estimating}. Note that as the vaccine efficacy increases, vaccine rate increases, or contact rate decreases, the effective reproduction number decreases.

\subsection{Basic COVID-19 Model} \label{sec:BC}
The next example was laid out as a basic COVID-19 model to aid in forecasting once the appropriate compartments for the disease had been identified \cite{gumel2021primer}. This example introduces more complex equations for $F$, involving additional interactions in $V$, and applies it to a specific disease. By applying our method to a compartmental COVID-19 model, we demonstrate its practical utility in analyzing contemporary and complex epidemiological challenges, where multiple infectious states with varying transmission dynamics play a crucial role. The compartments $S,E, I_a, I_s, I_h,$ and $R$ represent susceptible, exposed, infected asymptomatic, infected symptomatic, infected hospitalized, and recovered, respectively. The parameters $\beta_a, \beta_S,$ and $\beta_h$ represent the rates of infection due to asymptomatic, symptomatic, and hospitalized individuals, respectively. The parameters $\sigma, r, \gamma_a, \gamma_s, \varphi_S, \delta_s, \gamma_h,$ and $ \delta_h$ represent rate of becoming infectious, the percentage of asymptomatic people, rate of recovery for asymptomatic infected, rate of recovery for symptomatic infected, rate of hospitalization of symptomatic, death rate of symptomatic, rate of recovery for hospitalized and death rate of hospitalized, respectively. The Basic COVID-19 system is described by the following system of equations: 

\begin{align}
\frac{dS}{dt} &= -\frac{\beta_a I_a+\beta_S I_S+\beta_h I_h}{N}S,\label{eq:S_BC} \\
\frac{dE}{dt} &= \frac{\beta_a I_a+\beta_S I_S+\beta_h I_h}{N}S-\sigma E,\label{eq:E_BC}
\end{align}
\begin{align}
\frac{dI_a}{dt} &= r\sigma E-\gamma_a I_a,\label{eq:Ia_BC} \\
\frac{dI_S}{dt} &= (1-r)\sigma E-(\varphi_s+\gamma_s+\delta_s) I_s,\label{eq:Is_BC} \\
\frac{dI_h}{dt} &= \varphi_s I_s -(\gamma_h+\delta_h)I_h,\label{eq:Ih_BC} \\
\frac{dR}{dt} &= \gamma_a I_a+\gamma_s I_s+\gamma_h I_h. \label{eq:R_BC} 
\end{align}

\begin{figure}[H]
    \includegraphics[width=0.7\textwidth]{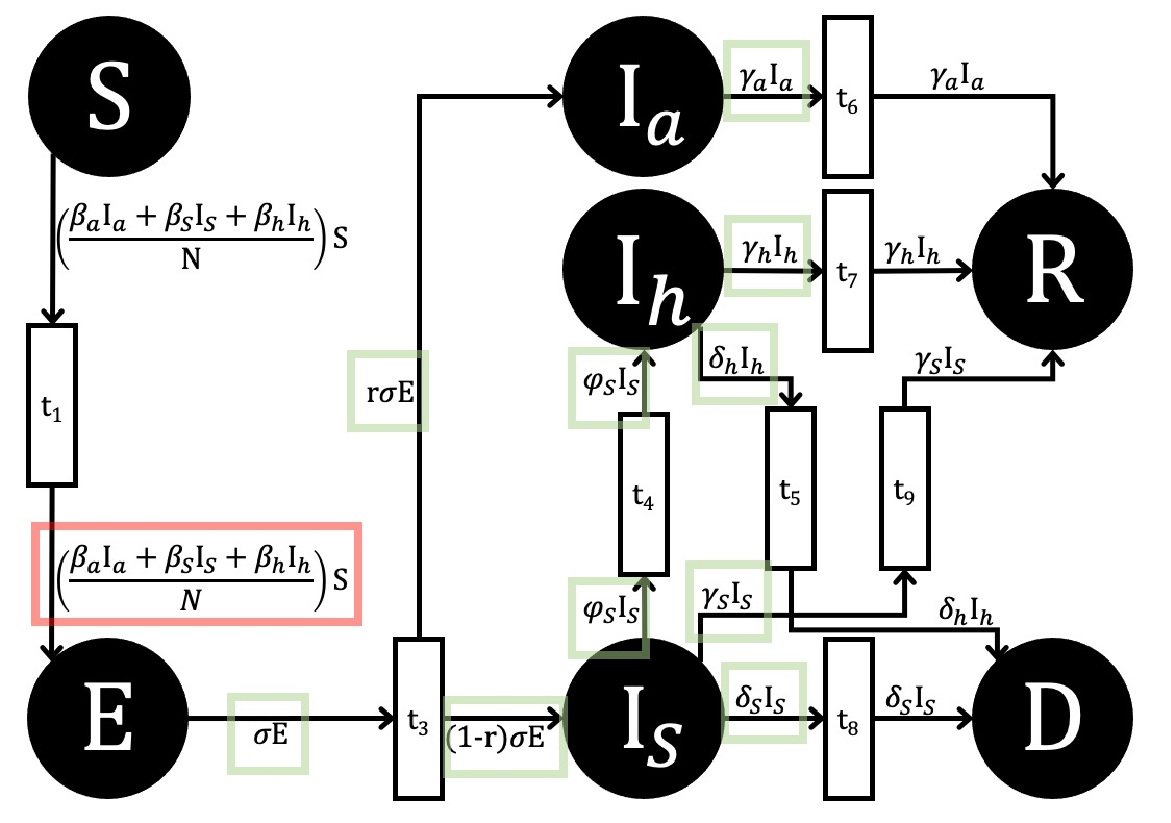}
    \caption{{A variable} 
 arc weight COVID-19 Petri Net mapped from Equations \eqref{eq:S_BC}--\eqref{eq:R_BC}. {The arc weights used in forming matrix $F_i(x)$ are outlined in red, and the arc weights used in forming matrix $V_i(x)$ are outlined in green.}}
    \label{fig:Basic Covid PN Lay}
\end{figure}

The Equations \eqref{eq:S_BC}--\eqref{eq:R_BC} are mapped to a variable arc weight Petri Net to yield Figure~\ref{fig:Basic Covid PN Lay}.

From Figure \ref{fig:Basic Covid PN Lay}, we can see that there are four infected places ($E, I_a, I_S, I_h$) which are the first, second, third, and fourth infected compartments, respectively. Our resulting $F$ and $V$ matrices will therefore be four by four. The disease free equilibrium is given by $(S^*, E^*, I_a^*, I_s^*, I_h^*, R^*)=(S^*, 0, 0, 0, R^*)$. We defined $N=S+E+I_a+I_S+I_h+R$, so $N^*=S^*+R^*$. By inspecting the PN model, we see that the term $\frac{\beta_a I_a+\beta_S I_S+\beta_h I_h}{N}S$ is the only way the Susceptible place leads to an infected place. The transmission matrix starts with the resulting $\mathcal{F}$ using Equation \eqref{eq:NGMPN_F_Transm} is
\begin{align*}
\mathcal{F}_{1,2,3,4}(x) &=\begin{bmatrix}
\frac{\beta_a I_a+\beta_S I_S+\beta_h I_h}{N}S & \frac{\beta_a I_a+\beta_S I_S+\beta_h I_h}{N}S & \frac{\beta_a I_a+\beta_S I_S+\beta_h I_h}{N}S & \frac{\beta_a I_a+\beta_S I_S+\beta_h I_h}{N}S \\
0 & 0 & 0 & 0 \\
0 & 0 & 0 & 0 \\
0 & 0 & 0 & 0 
\end{bmatrix} .
\end{align*}

Taking the partial derivative with respect to $E, I_a, I_S, I_h$ for the first, second, third, and fourth columns, respectively, gives us 
\begin{align*}
F &=\begin{bmatrix}
0 & \frac{\beta_a S^*}{N^*} & \frac{\beta_S S^*}{N^*} & \frac{\beta_h S^*}{N^*} \\
0 & 0 & 0 & 0 \\
0 & 0 & 0 & 0 \\
0 & 0 & 0 & 0 
\end{bmatrix} ,
\end{align*}
where $S^*, N^*$ are the $S,N$ values at the DFE.

\textls[-15]{For $V$, the transition matrix, we start by finding $\mathcal{V}$ based on the arc weights between infected places and the arc weights leaving the infected places as shown in Equation \eqref{eq:NGMPN_V}.
Thus,}
\begin{align*}
\mathcal{V}_{1,2,3,4}(x) &=\begin{bmatrix}
\sigma E & 0 & 0 & 0 \\
-r\sigma E & \gamma_a I_a & 0 & 0 \\
-(1-r)\sigma E & 0 & \varphi_S I_S+\gamma_S I_S +\delta_S I_S & 0 \\
0 & 0 & -\varphi_S I_S & \gamma_h I_h +\delta_h I_h
\end{bmatrix} .
\end{align*}
Then, taking the partial derivative with respect to $E, I_a, I_S, I_h$ for the first, second, third, and fourth columns, respectively, gives us 
\begin{align*}
V &=\begin{bmatrix}
\sigma & 0 & 0 & 0 \\
-r\sigma  & \gamma_a  & 0 & 0 \\
-(1-r)\sigma  & 0 & \varphi_S +\gamma_S +\delta_S& 0 \\
0 & 0 & \varphi_S & \gamma_h+\delta_h 
\end{bmatrix} .
\end{align*}

Hence,
\begin{align*}
    R_0=\varphi(FV^{-1})=\left(\left(\frac{\beta_a r}{\gamma_a}\right) + \left(\frac{\beta_s (1-r)}{\varphi_s+\gamma_s+\delta_s}\right) + \left(\frac{\beta_h (1-r)\varphi_s}{(\varphi_s+\gamma_s+\delta_s)(\gamma_h+\delta_h)}\right)\right)\left(\frac{S^*}{N^*}\right).
\end{align*}

The values for the Petri Net $F$ and $V$ are equivalent to those values for the ODE $F$ and $V$, and the resulting $R_0$ matches the $R_0$ for the ODE system in the Gumel paper \cite{gumel2021primer} seen in Equations (2) and (3).

\subsection{{Nonlinear System}
} \label{sec:NL}
The paper by Rohith et al. \cite{rohith2020dynamics} involves nonlinear compartment-to-compartment dynamics. These types of dynamics are more applicable when considering the social aspects of disease spread, drug pharmacokinetics, or disease cell spread within an individual. The analysis of the nonlinear system underscores the robustness of the NGMPN framework, demonstrating its applicability even when departing from standard mass-action assumptions to incorporate more realistic, behavior-driven transmission dynamics. The parameters $\mu, \beta, \alpha, \sigma,$ and $\gamma$ represent the birth/death rate, the per capita contact rate, the psychological or inhibitory effect, the rate of infectiousness, and the recovery rate, respectively. The nonlinear system is described by the following system of equations:
\begin{align}
\frac{dS}{dt} &= \mu - \frac{\beta S I}{1+\alpha I^2}-\mu S,\label{eq:NL_SEIR_S} \\
\frac{dE}{dt} &= \frac{\beta S I}{1+\alpha I^2}-\sigma E-\mu E,\label{eq:NL_SEIR_E} \\
\frac{dI}{dt} &= \sigma E -\gamma I -\mu I, \label{eq:NL_SEIR_I} \\
\frac{dR}{dt} &= \gamma I -\mu R.\label{eq:NL_SEIR_R} 
\end{align}

\begin{figure}[H]
    \includegraphics[width=0.6\textwidth]{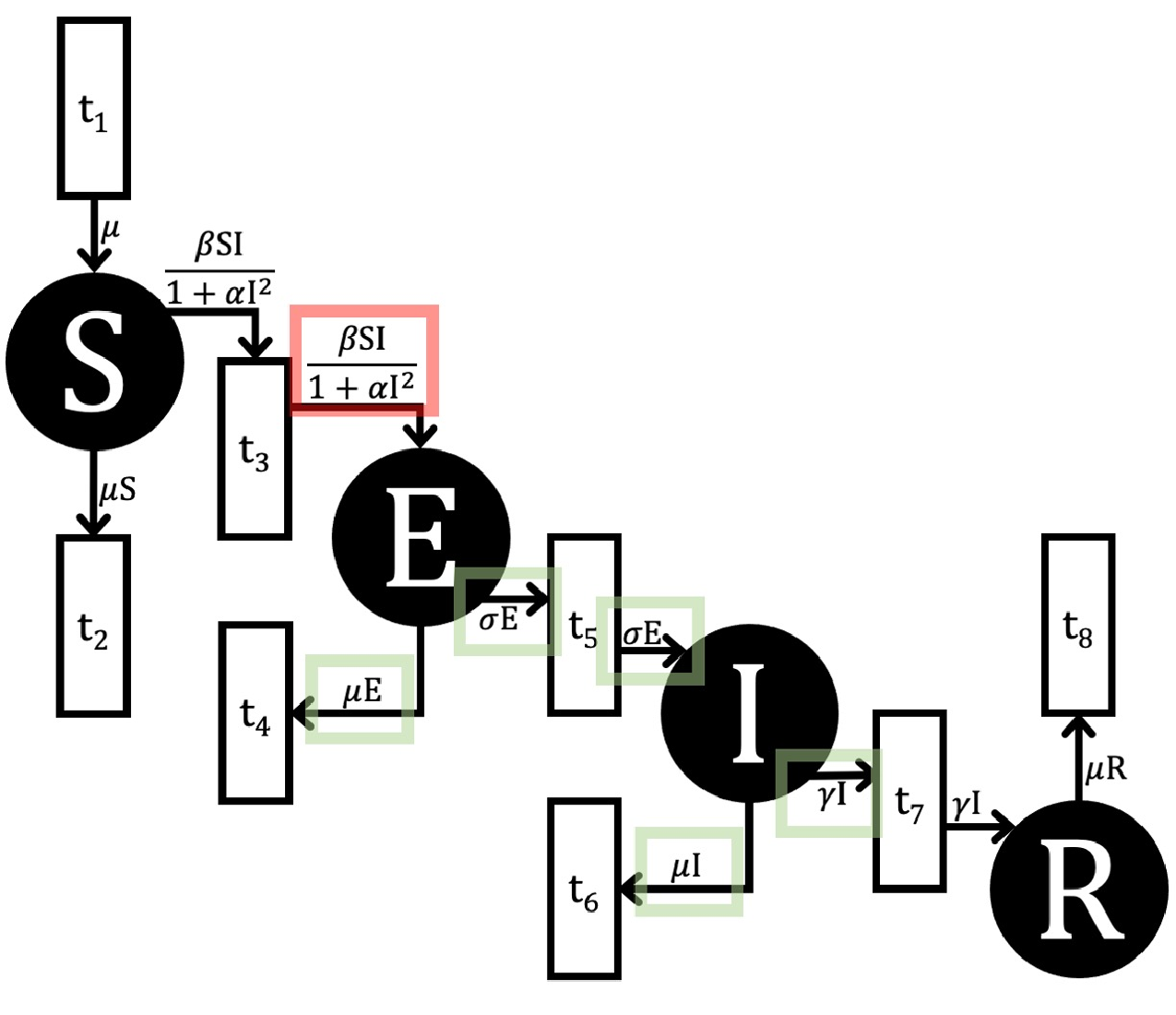}
    \caption{
 A variable arc weight nonlinear SEIR PN mapped from Equations (\ref{eq:NL_SEIR_S})--(\ref{eq:NL_SEIR_R}). {The arc weights used in forming matrix $F_i(x)$ are outlined in red, and the arc weights used in forming matrix $V_i(x)$ are outlined in green.}}
    \label{fig:NL SEIR PN Lay}
\end{figure}

From {Figure} 
 \ref{fig:NL SEIR PN Lay} we can see that there are two infected places $(E,I)$. Thus, our resulting $F$ and $V$ matrices will be two by two. Then, from inspecting the PN model, we see that $\frac{\beta SI}{1+\alpha I^2}$ is the only way that the new infected population enters. Thus the resulting $\mathcal{F}$ is
\begin{align*}
\mathcal{F}_{1,2}(x) &=\begin{bmatrix}
\frac{\beta S I}{1+\alpha I^2} & \frac{\beta S I}{1+\alpha I^2} \\
0 & 0 
\end{bmatrix}. 
\end{align*}

\vspace{-19pt}

Then, finding the DFE of the PN system to be $(S^*, E^*, I^*, R^*)=(1,0,0,0)$. We then take the partial derivative of each element of $\mathcal{F}_{1,2}$ with respect to $E, I$ for the first and second columns, respectively, to give the transmission matrix $F$,
\begin{align*}
F &=\begin{bmatrix}
0 & \beta \\
0 & 0 
\end{bmatrix}. 
\end{align*}

For the transition matrix $V$, we start by finding $\mathcal{V}$ based on the arc weights between infected places and the arc weights leaving the infected places, yielding
\begin{align*}
\mathcal{V}_{1,2}(x) &=\begin{bmatrix}
\sigma E + \mu E & 0 \\-\sigma E & \gamma I + \mu I
\end{bmatrix}. 
\end{align*}
Then, taking the partial derivative with respect to $E, I$ for the first and second column gives~us
\begin{align*}
V &=\begin{bmatrix}
\sigma+ \mu  & 0 \\
-\sigma  & \gamma  + \mu 
\end{bmatrix}. 
\end{align*}
Then we find the inverse of $V$,
\begin{align*}
V^{-1} &=\begin{bmatrix}
\frac{1}{\sigma+ \mu}  & 0 \\
\frac{\sigma}{(\sigma+ \mu)(\gamma  + \mu)}  & \frac{1}{\gamma  + \mu}  
\end{bmatrix}. 
\end{align*}

Then with the values of $F$ and $V^{-1}$ we can find the dominant eigen value of $(FV^{-1})$:
\begin{align*}
    R_0=\rho(FV^{-1})=\frac{\sigma \beta}{(\sigma+\mu)(\gamma+\mu)},
\end{align*}
which is the same $R_0$ value found in the Rohith paper \cite{rohith2020dynamics}. This nonlinear example, with its functional transmission term of  
$\frac{\beta SI}{1+\alpha I^2}$ represented by a single arc in Figure \ref{fig:NL SEIR PN Lay}, perfectly illustrates the strength of our framework in handling VAPN models. A geometric, path-summing method based on constant-rate transitions would be difficult to apply to such a compact representation.

\subsection{Patch System Overview} \label{sec:Patch}
The paper {``Multi-patch and multi-group epidemic models: a new framework''} 
by Bichara and Iggidr \cite{bichara2018multi} presents a multi-patch and multi-group model for infectious disease transmission. The model accommodates interactions occurring across an arbitrary number of patches and groups where the infection propagates. The compartments $S_i, E_i, I_i, R_i$ still represent susceptible, exposed, infected, and recovered, respectively, just now for patch $i$. Similarly, the parameters $\Pi, \beta, \mu, \eta, \upsilon, \gamma,$ and $ \delta $ represent birth/immigration, rate of infection, natural death rate, loss of immunity rate, rate of becoming infectious, rate of recovery, and rate of disease-related death, respectively. The parameter $\beta$ is a patch specific term, while $\Pi, \mu, \eta, \upsilon, \gamma,$ and $ \delta$  are group specific terms. People of group $i$ spend on average some time in patch $j$,$j=1,\ldots,n$. The susceptible, exposed, infected, and recovered populations of group $i$ spend $m_{ij}, n_{ij}, p_{ij},$ and $q_{ij}$ proportions of time in patch\linebreak  $j$ for $j=1,\ldots,n$.
The ODEs are outlined as
\begin{align}
\frac{dS_i}{dt} &= \Pi_i-\sum_{j=1}^{n} \beta_j m_{ij}S_i \frac{\sum_{k=1}^{u}p_{kj}I_k}{\sum_{k=1}^{u}(m_{kj}S_k+n_{kj}E_k+p_{kj}I_k+q_{kj}R_k)}-\mu_i S_i +\eta_i R_i,\label{eq:Si_Pa} \\
\frac{dE_i}{dt} &= \sum_{j=1}^{n} \beta_j m_{ij}S_i \frac{\sum_{k=1}^{u}p_{kj}I_k}{\sum_{k=1}^{u}(m_{kj}S_k+n_{kj}E_k+p_{kj}I_k+q_{kj}R_k)}-\upsilon_i E_I-\mu_i E_i,\label{eq:Ei_Pa} \\
\frac{dI_i}{dt} &= \upsilon_i E_I-\gamma_i I_i -\delta_i I_i-\mu_i I_i,\label{eq:Ii_Pa} \\
\frac{dR_i}{dt} &= \gamma_i I_i-\eta_i R_i-\mu_i R_i. \label{eq:Ri_Pa} 
\end{align}

\subsubsection*{{Two Patch System}
} 
We take the ODEs from Bichara and Iggidr \cite{bichara2018multi} Equations \eqref{eq:Si_Pa}--\eqref{eq:Ri_Pa} and lay them out explicitly for two patches and two groups. Although the ODE and corresponding PN model are only applicable to this small patch and group size, this model can be extended to any size, with the interactions between any two given patches being equivalent to the system laid out below. The multi-patch model highlights the scalability of our approach to spatially structured populations, providing a tool to investigate the impact of geographic factors and population mobility on disease spread. The Two Patch system is described by the following system of equations:  
\begin{adjustwidth}{-\extralength}{0cm}
\begin{align}
\frac{dS_1}{dt} &= \Pi_1-
\left( \beta_1 m_{11}S_1 \frac{p_{11}I_1+p_{21}I_2}{(m_{11}S_1+n_{11}E_1+p_{11}I_1+q_{11}R_1)+(m_{21}S_2+n_{21}E_2+p_{21}I_2+q_{21}R_2)} \right. +\nonumber\\
&\qquad \left.  \beta_2 m_{12}S_1 \frac{p_{12}I_1+p_{22}I_2}{(m_{12}S_1+n_{12}E_1+p_{12}I_1+q_{12}R_1)+(m_{22}S_2+n_{22}E_2+p_{22}I_2+q_{22}R_2)} \right)  -\mu_1 S_1 +\eta_1 R_1,\label{eq:S1_Pa} \\[1em]
\frac{dS_2}{dt} &= \Pi_2 -\left( \beta_1 m_{21}S_2 \frac{p_{11}I_1+p_{21}I_2}{(m_{11}S_1+n_{11}E_1+p_{11}I_1+q_{11}R_1)+(m_{21}S_2+n_{21}E_2+p_{21}I_2+q_{21}R_2)} \right. + \nonumber\\
&\qquad \left. \beta_2 m_{22}S_2 \frac{p_{12}I_1+p_{22}I_2}{(m_{12}S_1+n_{12}E_1+p_{12}I_1+q_{12}R_1)+(m_{22}S_2+n_{22}E_2+p_{22}I_2+q_{22}R_2)} \right)  -\mu_2 S_2 +\eta_2 R_2,\label{eq:S2_Pa} \\[1em]
\frac{dE_1}{dt} &= \left( \beta_1 m_{11}S_1 \frac{p_{11}I_1+p_{21}I_2}{(m_{11}S_1+n_{11}E_1+p_{11}I_1+q_{11}R_1)+(m_{21}S_2+n_{21}E_2+p_{21}I_2+q_{21}R_2)} \right. + \nonumber\\
&\qquad \left.  \beta_2 m_{12}S_1 \frac{p_{12}I_1+p_{22}I_2}{(m_{12}S_1+n_{12}E_1+p_{12}I_1+q_{12}R_1)+(m_{22}S_2+n_{22}E_2+p_{22}I_2+q_{22}R_2)} \right)  -\upsilon_1 E_1-\mu_1 E_1,\label{eq:E1_Pa}
\end{align}
\begin{align}
\frac{dE_2}{dt} &= \left( \beta_1 m_{21}S_2 \frac{p_{11}I_1+p_{21}I_2}{(m_{11}S_1+n_{11}E_1+p_{11}I_1+q_{11}R_1)+(m_{21}S_2+n_{21}E_2+p_{21}I_2+q_{21}R_2)} \right.+ \nonumber\\
&\qquad \left. \beta_2 m_{22}S_2 \frac{p_{12}I_1+p_{22}I_2}{(m_{12}S_1+n_{12}E_1+p_{12}I_1+q_{12}R_1)+(m_{22}S_2+n_{22}E_2+p_{22}I_2+q_{22}R_2)} \right) -\upsilon_2 E_2-\mu_2 E_2,\label{eq:E2_Pa} \\[1em]
\frac{dI_1}{dt} &= \upsilon_1 E_1-\gamma_1 I_1 -\delta_1 I_1-\mu_1 I_1,\label{eq:I1_Pa} \\
\frac{dI_2}{dt} &= \upsilon_2 E_2-\gamma_2 I_2 -\delta_2 I_2-\mu_2 I_2,\label{eq:I2_Pa} \\
\frac{dR_1}{dt} &= \gamma_1 I_1-\eta_1 R_1-\mu_1 R_1 ,\label{eq:R1_Pa}\\
\frac{dR_2}{dt} &= \gamma_2 I_2-\eta_2 R_2-\mu_2 R_2. \label{eq:R2_Pa}
\end{align}
\end{adjustwidth}

\begin{figure}[H]
    \includegraphics[width=0.7\textwidth]{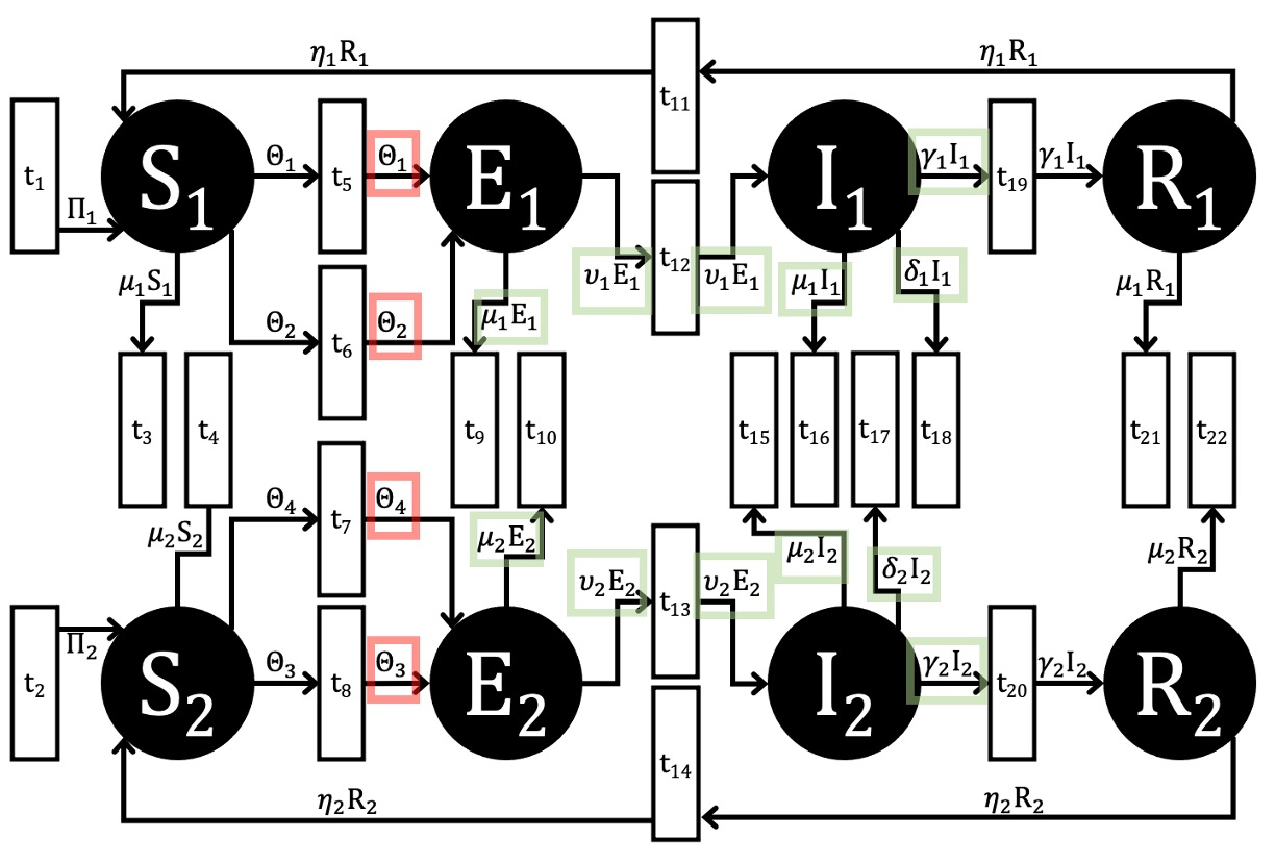}
    \caption{ 
 A variable arc weight two patch SEIR PN mapped from Equations (\ref{eq:S1_Pa})--(\ref{eq:R2_Pa}). {The arc weights used in forming matrix $F_i(x)$ are outlined in red, and the arc weights used in forming matrix $V_i(x)$ are outlined in green.}}
    \label{fig:Patch PN}
\end{figure}

From Figure \ref{fig:Patch PN} we find there are two infected places ($E$, $I$) per patch. Thus, with two patches, our resulting $F$ and $V$ matrices will be four by four. 
\vspace{-12pt}

{To simplify, we define}
\begin{align*}
\small 
\Theta_1=\beta_1 m_{11}S_1 \frac{p_{11}I_1+p_{21}I_2}{(m_{11}S_1+n_{11}E_1+p_{11}I_1+q_{11}R_1)+(m_{21}S_2+n_{21}E_2+p_{21}I_2+q_{21}R_2)},\\
\Theta_2=\beta_2 m_{12}S_1 \frac{p_{12}I_1+p_{22}I_2}{(m_{12}S_1+n_{12}E_1+p_{12}I_1+q_{12}R_1)+(m_{22}S_2+n_{22}E_2+p_{22}I_2+q_{22}R_2)},\\
\Theta_3=\beta_1 m_{21}S_2 \frac{p_{11}I_1+p_{21}I_2}{(m_{11}S_1+n_{11}E_1+p_{11}I_1+q_{11}R_1)+(m_{21}S_2+n_{21}E_2+p_{21}I_2+q_{21}R_2)},\\
\Theta_4=\beta_2 m_{22}S_2 \frac{p_{12}I_1+p_{22}I_2}{(m_{12}S_1+n_{12}E_1+p_{12}I_1+q_{12}R_1)+(m_{22}S_2+n_{22}E_2+p_{22}I_2+q_{22}R_2)}.
\end{align*}

Then we obtain the DFE by setting the infected places token levels to zero, setting the net arc weight of non-infected places to 0, and solving for the non-infected places token values $(S_1^*, S_2^*, E_1^*, E_2^*, I_1^*, I_2^*, R_1^*, R_2^*)=(\frac{\Pi_1}{\mu_1},\frac{\Pi_2}{\mu_2}, 0, 0, 0, 0, 0, 0)$. Then the matrix $\mathcal{F}$ by looking at the arc weights of arcs going into the infected compartments, giving us 
 \begin{align*}
\mathcal{F}_{1,2,3,4}(x) &=\begin{bmatrix}
\Theta_1+\Theta_2 & \Theta_1+\Theta_2 & \Theta_1+\Theta_2 & \Theta_1+\Theta_2\\
\Theta_3+\Theta_4 & \Theta_3+\Theta_4 & \Theta_3+\Theta_4 & \Theta_3+\Theta_4\\
0 & 0 & 0 & 0\\
0 & 0 & 0 & 0
\end{bmatrix}. 
\end{align*}
\textls[--25]{Finding the rate of change with respect to the infected states and evaluating at the DFE~yields }
\begin{align*}
F &=\begin{bmatrix}
0 & 0 & \frac{\beta_1 m_{11}p_{11}S_1}{m_{11}S_1+m_{21}S_2}+\frac{\beta_2 m_{12}p_{12}S_1}{m_{12}S_1+m_{22}S_2} & \frac{\beta_1 m_{11}p_{12}S_1}{m_{11}S_1+m_{21}S_2}+\frac{\beta_2 m_{12}p_{22}S_1}{m_{12}S_1+m_{22}S_2} \vspace{8pt}\\
0 & 0 & \frac{\beta_1 m_{21}p_{11}S_2}{m_{11}S_1+m_{21}S_2}+\frac{\beta_2 m_{22}p_{12}S_2}{m_{12}S_1+m_{22}S_2} & \frac{\beta_1 m_{21}p_{21}S_2}{m_{11}S_1+m_{21}S_2}+\frac{\beta_2 m_{22}p_{22}S_2}{m_{12}S_1+m_{22}S_2}\\
0 & 0 & 0 & 0\\
0 & 0 & 0 & 0
\end{bmatrix}.
\end{align*}
{For $\mathcal{V}$, we find arcs leaving the infected places and transitioning between infected places,~yielding }
\begin{align*}
\mathcal{V}_{1,2,3,4}(x) &= \begin{bmatrix}
\upsilon_1 E_1 + \mu_1 E_1 & 0 & 0 & 0\\
0 & \upsilon_2 E_2 + \mu_2 E_2  & 0 & 0\\
-\upsilon_1 E_1 & 0 & \gamma_1 I_1 + \delta_1 I_1 + \mu_1 I_1 & 0\\
0 & -\upsilon_2 E_2 & 0 & \gamma_2 I_2 + \delta_2 I_2 + \mu_2 I_2
\end{bmatrix}.
\end{align*}
We then find the rate of change with respect to the infected place $I$, and evaluate at the DFE. From this, we find the inverse to give
\begin{align*}
V &= \begin{bmatrix}
\upsilon_1 + \mu_1 & 0 & 0 & 0\\
0 & \upsilon_2 + \mu_2  & 0 & 0\\
-\upsilon_1 & 0 & \gamma_1 + \delta_1 + \mu_1 & 0\\
0 & -\upsilon_2 & 0 & \gamma_2 + \delta_2 + \mu_2
\end{bmatrix}, \\
V^{-1} &= \begin{bmatrix}
\frac{1}{\upsilon_1 + \mu_1} & 0 & 0 & 0\\
0 & \frac{1}{\upsilon_2 + \mu_2} & 0 & 0\\
\frac{\upsilon_1}{(\upsilon_1 + \mu_1)(\gamma_1 + \delta_1 + \mu_1)} & 0 & \frac{1}{\gamma_1 + \delta_1 + \mu_1} & 0\\
0 & \frac{\upsilon_2}{(\upsilon_2 + \mu_2)(\gamma_2 + \delta_2 + \mu_2)} & 0 & \frac{1}{\gamma_2 + \delta_2 + \mu_2}
\end{bmatrix}.
\end{align*}
Now that we have $F$ and $V$ matrices, we can find $R_0$ by multiplying and determining the dominant eigenvalue. Thus,
\vspace{-10pt}
\begin{adjustwidth}{-\extralength}{0cm}
\begin{align*}
R_0 &= \varphi\begin{bmatrix}
F \cdot V^{-1} 
\end{bmatrix}=  \varphi
{\small\begin{bmatrix}
\frac{\left(\frac{\beta_{1} m_{{11}} p_{{11}} S_{1}}{m_{{11}} S_{1}+m_{{21}} S_{2}}+\frac{\beta_{2} m_{{12}} p_{{12}} S_{1}}{S_{1} m_{{12}}+S_{2} m_{{22}}}\right) \upsilon_{1}}{\left(\upsilon_{1}+\mu_{1}\right) \left(\gamma_{1}+\delta_{1}+\mu_{1}\right)} & \frac{\left(\frac{\beta_{1} m_{{11}} p_{{12}} S_{1}}{m_{{11}} S_{1}+m_{{21}} S_{2}}+\frac{\beta_{2} m_{{12}} p_{{22}} S_{1}}{S_{1} m_{{12}}+S_{2} m_{{22}}}\right) \upsilon_{2}}{\left(\upsilon_{2}+\mu_{2}\right) \left(\gamma_{2}+\delta_{2}+\mu_{2}\right)} & \frac{\frac{\beta_{1} m_{{11}} p_{{11}} S_{1}}{m_{{11}} S_{1}+m_{{21}} S_{2}}+\frac{\beta_{2} m_{{12}} p_{{12}} S_{1}}{S_{1} m_{{12}}+S_{2} m_{{22}}}}{\gamma_{1}+\delta_{1}+\mu_{1}} & \frac{\frac{\beta_{1} m_{{11}} p_{{12}} S_{1}}{m_{{11}} S_{1}+m_{{21}} S_{2}}+\frac{\beta_{2} m_{{12}} p_{{22}} S_{1}}{S_{1} m_{{12}}+S_{2} m_{{22}}}}{\gamma_{2}+\delta_{2}+\mu_{2}} 
\vspace{8pt}\\
 \frac{\left(\frac{\beta_{1} m_{{21}} p_{{11}} S_{2}}{m_{{11}} S_{1}+m_{{21}} S_{2}}+\frac{\beta_{2} m_{{22}} p_{{12}} S_{2}}{S_{1} m_{{12}}+S_{2} m_{{22}}}\right) \upsilon_{1}}{\left(\upsilon_{1}+\mu_{1}\right) \left(\gamma_{1}+\delta_{1}+\mu_{1}\right)} & \frac{\left(\frac{\beta_{1} m_{{21}} p_{{21}} S_{2}}{m_{{11}} S_{1}+m_{{21}} S_{2}}+\frac{\beta_{2} m_{{22}} p_{{22}} S_{2}}{S_{1} m_{{12}}+S_{2} m_{{22}}}\right) \upsilon_{2}}{\left(\upsilon_{2}+\mu_{2}\right) \left(\gamma_{2}+\delta_{2}+\mu_{2}\right)} & \frac{\frac{\beta_{1} m_{{21}} p_{{11}} S_{2}}{m_{{11}} S_{1}+m_{{21}} S_{2}}+\frac{\beta_{2} m_{{22}} p_{{12}} S_{2}}{S_{1} m_{{12}}+S_{2} m_{{22}}}}{\gamma_{1}+\delta_{1}+\mu_{1}} & \frac{\frac{\beta_{1} m_{{21}} p_{{21}} S_{2}}{m_{{11}} S_{1}+m_{{21}} S_{2}}+\frac{\beta_{2} m_{{22}} p_{{22}} S_{2}}{S_{1} m_{{12}}+S_{2} m_{{22}}}}{\gamma_{2}+\delta_{2}+\mu_{2}} 
\\
 0 & 0 & 0 & 0 
\\
 0 & 0 & 0 & 0 
\end{bmatrix}}.
\nonumber
\end{align*}
\end{adjustwidth}
The dominant eigenvalue for this particular system is an incredibly long algebraic expression. However, with the $F$ and $V$ matrices matching that of the system outlined in the paper~\cite{bichara2018multi} and subsequently for two patches and two groups in Equations \eqref{eq:S1_Pa}--\eqref{eq:R2_Pa}, the resulting $R_0$ coming from the dominant eigenvalue will be equal.

\subsection{SIR Vector-Borne Model} \label{sec:VB}

The outline of the original ODE system, as described by Wedajo et al. \cite{wedajo2018analysis}, is based on the SIR Malaria model. However, the structure of this model can be taken more broadly to represent a simple example of vector-borne disease modeling. The SIR vector-borne model demonstrates the adaptability of our method to different modes of transmission, extending its applicability from directly transmitted diseases to complex vector-host systems. The compartments $ S_h$, $ I_h$, $ R_h$, $ S_v$, and $I_v$ represent susceptible humans, infected humans, recovered humans, susceptible vectors, and infected vectors, respectively. Note that in a vector-borne model, the vector's infected compartments are considered general infected compartments for purposes of calculating $R_0$. The vector for Malaria in the original model is mosquitoes. The parameters $\Pi, \beta_{hv}, \mu_h, \delta, \alpha, \sigma, \Lambda, \beta_{vh},$ and $ \mu_v$ represent birth/immigration of humans, the rate of infected vectors infecting humans, human natural death rate, the rate of infected humans immigrating, the rate of infected humans emigrating, the recovery rate, the birth/immigration rate of the vector, the rate of infected humans infecting vectors, and natural vector death rate, respectively. The SIR Vector-borne system is described by the following system of equations:  
\begin{align}
\frac{dS_h}{dt}&=\Pi -\beta_{hv}S_h I_v- \mu_h S_h, \label{eq:Sh_SIR_Mal} \\
\frac{dI_h}{dt}&=\beta_{hv}S_h I_v +\delta I_h \ldots \nonumber\\&-(\alpha + \sigma + \mu_h) I_h, \label{eq:Ih_SIR_Mal} \\
\frac{dR_h}{dt}&=\sigma I_h - \mu_h R_h, \label{eq:Rh_SIR_Mal} \\
\frac{dS_v}{dt}&=\Lambda -\beta_{vh}S_v I_h-\mu_v S_v ,\label{eq:Sv_SIR_Mal} \\
\frac{dI_v}{dt}&=\beta_{vh}S_v I_h -\mu_v I_v. \label{eq:Iv_SIR_Mal} 
\end{align}

\begin{figure}[H]
    \includegraphics[width=0.7\textwidth]{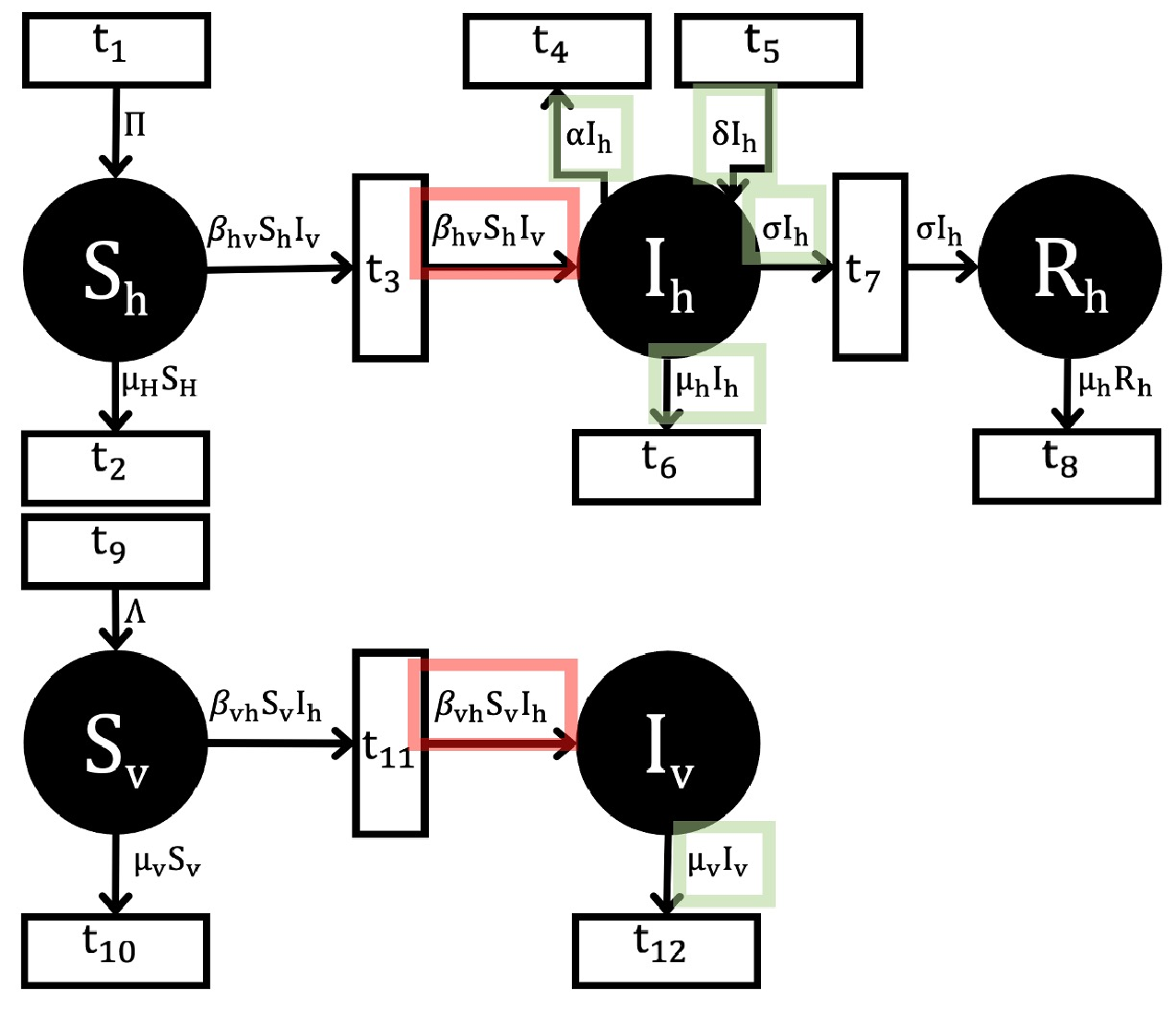}
    \caption{
A variable arc weight SIR Vector-borne PN mapped from Equations (\ref{eq:Sh_SIR_Mal})--(\ref{eq:Iv_SIR_Mal}). The arc weights used in forming matrix $F_i(x)$ are outlined in red, and the arc weights used in forming matrix $V_i(x)$ are outlined in green. }
    \label{fig:Vector_borne PN}
\end{figure}

From Figure \ref{fig:Vector_borne PN} we determine that there are two infected places ($I_H$, $I_v$) per patch. Then, $F$ and $V$ matrices will be two by two.

First, we obtain the DFE by setting the infected places token levels to zero, setting the net arc weight of non-infected places to 0, and solving for the non-infected places token values, $(S_h^*, S_v^*, I_h^*, I_v^*, R_h^*)=(\frac{\Pi}{\mu_h}, \frac{\Lambda}{\mu_v}, 0, 0, 0)$. Then the matrix $\mathcal{F}$ is obtained by looking at the arc weights of arcs going into the infected compartments, giving us 
 \begin{align*}
\mathcal{F}_{1,2}(x) &=\begin{bmatrix}
\beta_{hv}S_h I_v & \beta_{hv}S_h I_v \\
\beta_{vh}S_v I_h & \beta_{vh}S_v I_h
\end{bmatrix}. 
\end{align*}
Then we find the rate of change with respect to $I_h$ in the first column and $I_v$ in the second column, and plugging in the DFE to give
\begin{align*}
F &=\begin{bmatrix}
0 & \beta_{hv}\frac{\Pi}{\mu_h} \\
\beta_{vh}\frac{\Lambda}{\mu_v} & 0 
\end{bmatrix}.
\end{align*}
For $\mathcal{V}$, we look at arc weights of arcs leaving the infected places and transitioning from one infected place to another, yielding
\begin{align*}
\mathcal{V}_{1,2}(x) &= \begin{bmatrix}
\alpha I_h +\mu_h I_h +\sigma I_h -\delta I_h & 0 \\
0 & \mu_v I_v 
\end{bmatrix}.
\end{align*}
We then find the rate of change with respect to the infected place $I$, and evaluate at the DFE. From this, we find the inverse to give
\begin{align*}
V &= \begin{bmatrix}
\alpha +\mu_h +\sigma -\delta  & 0 \\
0 & \mu_v
\end{bmatrix}, \\
V^{-1} &= \begin{bmatrix}
\frac{1}{\alpha +\mu_h +\sigma -\delta}  & 0 \\
0 & \frac{1}{\mu_v}
\end{bmatrix}.
\end{align*}
Then the $F$ and $V$ matrices, allow us to determine $R_0$. Thus,
\begin{align*}
R_0 &= \varphi\begin{bmatrix}
F \cdot V^{-1} 
\end{bmatrix} = \sqrt{(\frac{\beta_{hv}\Pi}{\mu_h \mu_v}\cdot\frac{\beta_{vh}\Lambda}{(\mu_v)(\alpha +\mu_h +\sigma -\delta)})}.
\end{align*}

\textls[-20]{Some of the parameters are renamed in the layout of the ODE model in Equations~\eqref{eq:Sh_SIR_Mal}--\eqref{eq:Iv_SIR_Mal} versus the original \cite{wedajo2018analysis}. With this parameter renaming, the resulting $R_0$ value equals the exact $R_0$ expression in the Wedajo paper. }

\section{Numerical Verification} \label{sec:Numerical_Verification}
This section presents the numerical verification of the algebraic $R_0$ expressions derived from our framework. This validation is performed by comparing the analytical $R_0$, found using the NGMPN method, to a numerically estimated $R_0$ derived from simulation output. We specifically focus on simulation output from the deterministic Variable Arc Weight Petri Net (VAPN) formalism. The VAPN's ability to be numerically convergent with standard ODEs has been systematically studied, establishing it as a valid, deterministic PN implementation for this type of comparison \cite{reckell_numerical_2024}. Therefore, our numerical verification scheme is a direct, deterministic comparison, for which RRMSE is the most appropriate measure of error. While our NGMPN framework also applies to stochastic SPN models, numerical validation with SPN models would involve stochastic analysis, which would be suitable for future work. The numerical results presented offer an independent validation of our $R_0$ expression.


To achieve this validation, we conducted simulations of the Petri Net model using GPenSim, a general-purpose Petri Net simulator. GPenSim allows for close representation of ODE systems in Petri Nets \cite{reckell_numerical_2024}. By implementing the VAPN and analyzing the output of these simulations, specifically the progression of the susceptible and infectious populations under various parameter sets, we can numerically estimate the $R_0$ exhibited by the simulated system. This analysis is independent of any algebraic $R_0$ methods that look at the parent system. The subsequent results will demonstrate the concordance between the $R_0$ value obtained through algebraic methods and the $R_0$ value inferred from the GPenSim simulation outcomes, thereby verifying the accuracy and reliability of our analytical expression. To illustrate this, we compare the $R_0$ values obtained algebraically from the systems with those derived through various methods of analyzing the resulting simulation data. The comparison is made via the relative root mean square error (RRMSE).

\begin{figure} [H]
   \hspace{17mm}  \textbf{Petri Net $R_0$ vs Algebraic $R_0$ for SIRS system}\par
  \hspace{-7mm}  \includegraphics[width=0.8\linewidth]{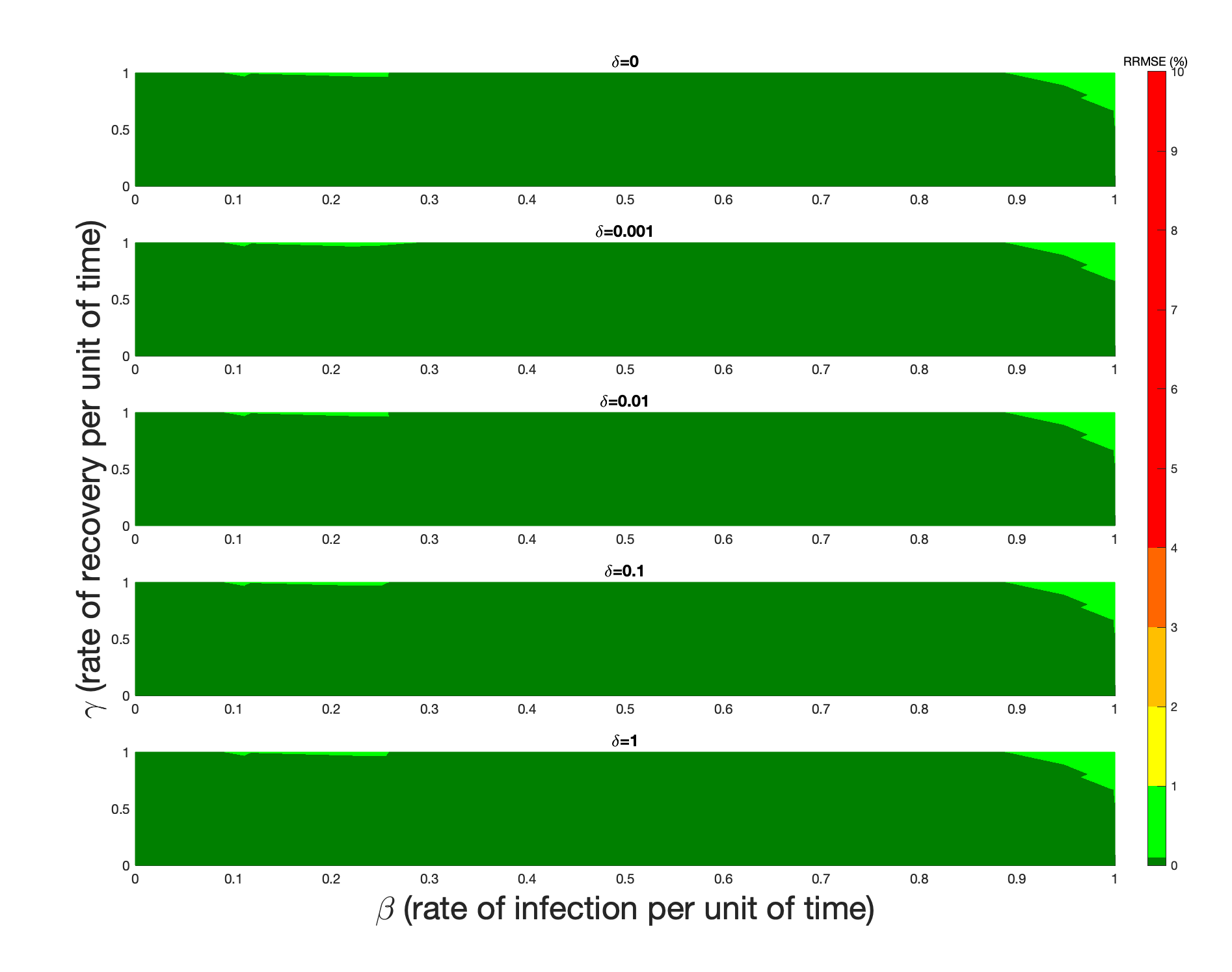}
    \caption{{Using} 
 the R library {``R0: Estimation of R0 and Real-Time Reproduction Number from Epidemics''} from Boelle and Obadia \cite{obadia2012r0, Boelle_Obadia} with estimation method of $R_0$ of Attack rate. The error of the simulations remains below 1\% RRMSE for all parameter values in the SIRS system and below 0.1\% RRMSE for most parameter sets and below 1\% RRMSE for all parameter sets, including all biologically plausible values. }
\label{fig:SIRS_PN_R0}
\end{figure}

\begin{figure}[H]
    \hspace{13mm} \textbf{Petri Net $R_0$ vs Algebraic $R_0$ for Nonlinear system}\par
  \hspace{-7mm}  \includegraphics[width=0.83\linewidth]{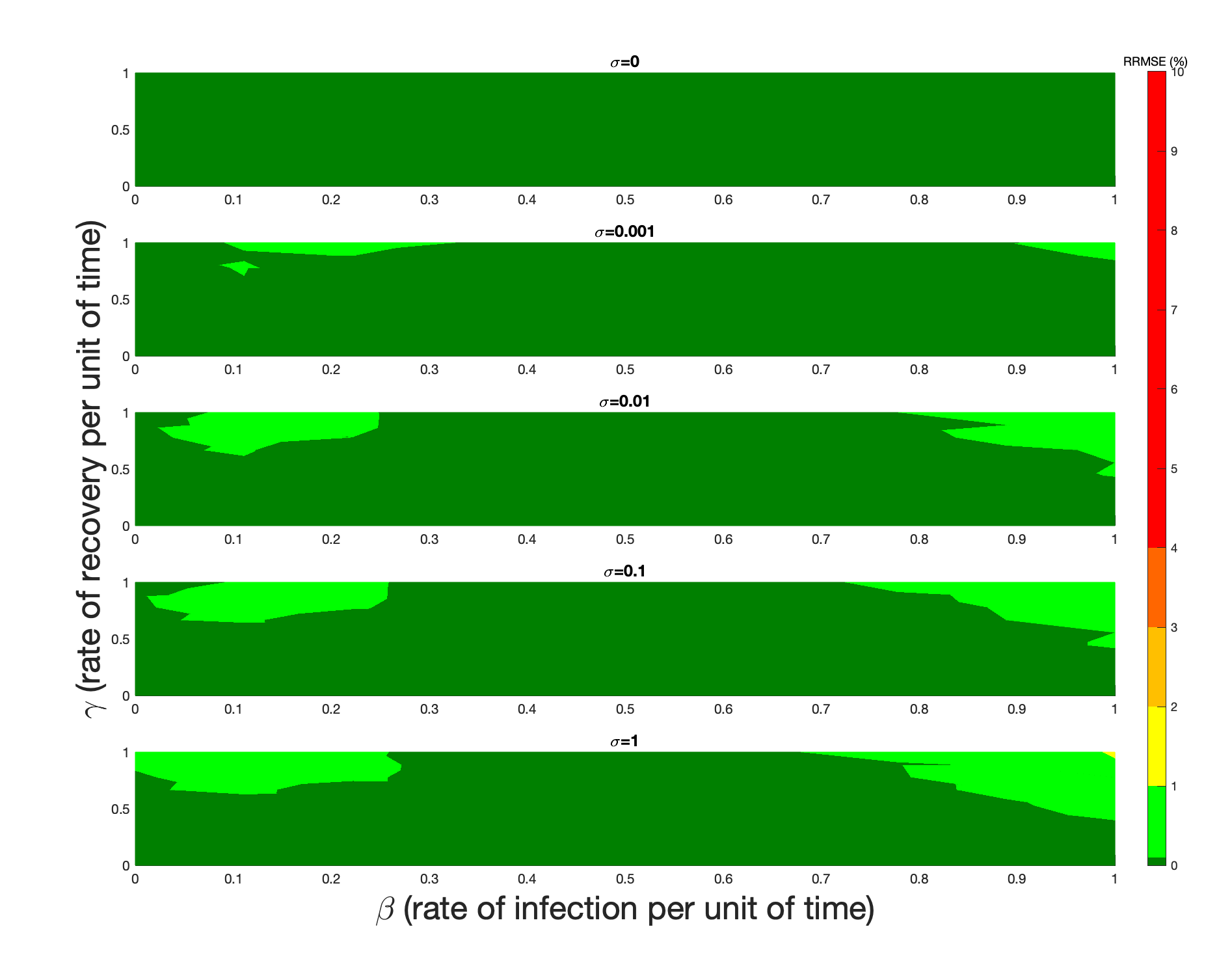}
    \caption{{Using} 
 the R library {{R0: Estimation of R0 and Real-Time Reproduction Number from Epidemics}} from Boelle and Obadia \cite{obadia2012r0, Boelle_Obadia} with estimation method of $R_0$ of ML. The error of the simulations remains below 1.2\% RRMSE for all parameter values in the Nonlinear system and below 0.1\% RRMSE for nearly all parameter values, including all biologically plausible values. }
\label{fig:NL_PN_R0}
\end{figure}

The Petri Net models and ODE models have been defined and run for a range of parameters. For models SIRS {Section} 
 \ref{sec:SIRS} and the Nonlinear model {Section} 
 \ref{sec:NL}, error surfaces were found for the $R_0$ values depending on three of the parameters.
The $R_0$ values from the VAPN system are found directly from the resulting simulation data using the {R package} 
{``R0: Estimation of R0 and Real-Time Reproduction Number from Epidemics'' Version:1.3-1} 
 library from Boelle and Obadia \cite{obadia2012r0, Boelle_Obadia}. The Maximum Likelihood (ML)  \cite{white2009estimation} method within the package was utilized to estimate $R_0$ from the simulation time series. This estimated $R_0$ is then compared to the analytically derived $R_0$ value from our NGMPN method. The error is quantified using RRMSE. We find strong agreement in numerical terms because the Petri Net $R_0$ calculations achieve RRMSE values consistently below 1--2\% for biologically plausible parameter sets, as shown in Figures \ref{fig:SIRS_PN_R0} and \ref{fig:NL_PN_R0}. While the R packages \cite{obadia2012r0, Boelle_Obadia} can produce confidence intervals (CIs) for stochastic data, our primary validation here is a deterministic-to-deterministic comparison, making RRMSE the most direct and rigorous metric for goodness of fit.

The algebraic expression using the NGM and NGMPN for the SIRS {Section} 
 \ref{sec:SIRS} systems is $R_0=\frac{\beta N}{\gamma}$. The following Figure \ref{fig:SIRS_PN_R0} compares the algebraic expression with that found from applying the attack rate method within the R library to the resulting simulation data across parameter values. For the Nonlinear system in Section \ref{sec:NL}, the algebraic expression of $R_0$ using the NGM and NGMPN is $R_0=\rho(FV^{-1})=\frac{\sigma \beta}{(\sigma+\mu)(\gamma+\mu)}$. Since $R_0$ includes the parameters $\beta$, $ \gamma$, and $\sigma$, these will be used as the parameters for comparing $R_0$ values across a range for each.


%

As shown in Figure \ref{fig:NL_PN_R0}, the Petri Net implementation of the nonlinear model produces $R_0$ values that closely match the algebraic solution, with RRMSE consistently below 0.1\% for most parameter sets and below 1\% for all parameter sets. This high level of agreement demonstrates the accuracy and numerical stability of the NGMPN method, even in the presence of nonlinear transmission and recovery terms. These results support the validity of the approach for more complex and mechanistically detailed epidemiological models.

%


\section{Discussion and Conclusions}
In this work, we devised a formal framework for finding the basic reproduction number $R_0$ within systems modeled by Petri Nets. This framework draws upon well-established methodologies employed in the analysis of ODE models. The relevance of applying $R_0$ in the context of Petri Nets is particularly pronounced, given the increasing utilization of Petri Nets in epidemiological modeling to depict more intricate and dynamic scenarios of disease spread. These scenarios encompass multi-compartment models, spatially distributed populations, and frameworks that incorporate social and nonlinear dynamics. Our approach builds on the next-generation matrix method, which is well understood within the realm of ODEs. Our results extend the definition of $R_0$ to Petri Nets by interpreting token transitions and arc weights in a manner that aligns with the underlying epidemiological processes. 
As seen in Section \ref{subsec:NGMPN}, by focusing on a formal, algorithmic identification of arcs corresponding to new infections, $F$, contrasted with other transitional changes, $V$, our framework offers a direct path to calculating $R_0$. This formal clarity enables the 'by-inspection' application and color-coding in our figures, which provides a direct visual representation of our method's partition logic.

The primary contributions of our paper are threefold. First, we provide a clear and formalism--agnostic operational framework for constructing the NGM. It applies directly to Petri Nets regardless of their structure, and we demonstrate its effectiveness on both traditional SPNs and the more flexible VAPNs. Second, we apply this method to a wide range of complex, contemporary models, including multi-strain COVID-19 and multi-patch systems. Lastly, we offer the first, to our knowledge, numerical validation of an NGMPN method, showing that our analytical results align precisely with $R_0$ values estimated from simulation data. This validation is a crucial step in ensuring that these theoretical tools are reliable for practical application in epidemiological modeling.

In parallel to the work performed by \cite{segovia_petri_2025}, we demonstrate that the NGM methodology can be effectively employed in the context of Petri Nets by first defining relevant matrices that represent the rate of infection, $F$, and the rate of removal or transition out of compartments, $V$. Subsequently, we compute the dominant eigenvalue of the product of $FV^{-1}$. This approach is consistent with established NGM techniques for ODE frameworks, resulting in a robust mechanism for quantifying the basic reproduction number when applied to Petri Nets. The examples we laid out, including models such as SIRS, SVEIR, SEEIR, Basic COVID-19, Nonlinear, Patch, and SIR Vector-borne, serve to illustrate the versatility of this methodology across a spectrum of disease dynamics characterized by various place to transition structures and arc weights.

A primary challenge encountered in adapting the NGM framework to Petri Nets pertains to effectively translating the compartmental structures and transition rates—typically articulated through arc weights and token movement—into a form compatible with the mathematical formulations employed in the NGM. Nevertheless, once the transition and removal rates are contextualized in terms of Petri Net parameters using Diekmann's approach \cite{diekmann2010construction}, the ensuing process of calculating $R_0$ becomes straightforward and closely parallels that utilized for ODEs. This provides a robust and generalizable method for scrutinizing disease dynamics in models based on Petri Nets.

Moreover, the framework afforded by Petri Nets enables the integration of more complex features that are often challenging to model using traditional ODE approaches. These features include variable parameters, discrete events, and nonlinear transitions~\cite{peleg2005using}. The SVEIR model in Section \ref{sec:SVEIR}  demonstrates how this framework naturally extends beyond the static $R_0$ to calculate the time-dependent effective reproduction number $R_e$, incorporating features like vaccination rates \cite{van2002reproduction, cori2013new}. Such flexibility renders Petri Nets a compelling tool for accurately representing real-world epidemiological systems that exhibit complex interactions and dependencies, including co-infections, varying transmission rates, and spatial heterogeneity in disease dissemination.

The proposed method maintains practical limitations that also limit other methods of $R_0$, including those used for ODEs. These practical challenges persist in their application to extensive or highly intricate Petri Nets. Specifically, the computational demands associated with calculating the dominant eigenvalue of large matrices can be substantial, particularly when addressing large-scale disease models that encapsulate numerous compartments. Again, these problems currently exist with all methods of algebraic $R_0$ formulation, but they should be noted. Additionally, considerations pertaining to the discrete nature of token movement in Petri Nets may necessitate careful deliberation when modeling diseases characterized by continuous dynamics. 
For SPNs in particular, the deterministic $R_0=1$ threshold has known explanatory limitations, as it does not fully capture the probability of stochastic fadeout where an outbreak dies out by chance even if $R_0>1$, especially with small initial numbers of infected individuals. This ``gap'' between the deterministic threshold and stochastic behavior can be bridged by considering the system's scale. As demonstrated numerically in \cite{reckell_numerical_2024}, SPN models show numerical convergence to their corresponding deterministic ODE and VAPN counterparts as the total population size is rescaled to larger values. While not a theoretical proof, this numerical convergence indicates that the deterministic $R_0$ derived from our framework is at minimum a practical approximation for the limit for the stochastic process. Therefore, while our $R_0$ value does not predict stochastic fadeout in small populations, it remains the crucial and accurate metric for predicting the mean, large-scale outbreak potential of the stochastic system. A more systematic analysis of the probabilistic behavior around this $R_0=1$ threshold for different population scales is an important area for future investigation.


Future work may be directed towards enhancing computational techniques to accommodate larger and more intricate systems, as well as exploring the application of this method to empirical epidemiological data. Overall, this work contributes to bridging the divide between traditional ODE-based modeling approaches and Petri Net frameworks, thereby providing a more comprehensive toolkit for understanding and predicting disease propagation across various {contexts}
.

\vspace{6pt}

\authorcontributions{{Conceptualization: T.R., B.S. and P.J.; Methodology:  T.R.; Software,  T.R.; Validation:  T.R.; Formal Analysis:  T.R.; Investigation:  T.R.; Resources:  T.R.; Data Curation:  T.R.; Writing---original draft preparation:  T.R.; Writing---review and editing,  T.R., B.S. and P.J.; Visualization:  T.R.; Supervision, B.S. and P.J.; Project Administration: B.S. and P.J.; Funding Acquisition: B.S. and P.J. All authors have read and agreed to the published version of the manuscript.}
}

\funding{{Trevor} 
 Reckell is partially supported by NIH grant DMS-1615879. Petar Jevti\'{c} and Beckett Sterner are supported by NIH grant 5R01GM131405-02.}

 \institutionalreview{{Not applicable.} 
}

\informedconsent{{Not applicable.} 
}

\dataavailability{Data and code used to produce data used in this paper is available in the publicly accessible repository GitHub at \url{https://github.com/trevorreckell/NGMPN}  (accessed on 10 October 2025).
}

\acknowledgments{
{Generative AI was used solely to resolve technical R and LaTeX errors, refine the discussion section, and identify specific introduction citations. Generative AI was not used in any of the methodology, scientific content, or analysis. 
}}

\conflictsofinterest{The authors declare no conflicts of interest.}

\begin{adjustwidth}{-\extralength}{0cm}
\reftitle{References}

\PublishersNote{}
\end{adjustwidth}
\end{document}